\documentclass[journal,symbol*]{new-aiaa}
\usepackage[utf8]{inputenc}
\usepackage{textcomp}

\usepackage{graphicx}
\usepackage{amsmath}
\usepackage{mathtools}
\usepackage{subcaption}
\usepackage[version=4]{mhchem}
\usepackage{siunitx}
\usepackage{longtable,tabularx}
\usepackage{amssymb}

\usepackage{amsthm}
\theoremstyle{definition}

\newtheorem{assumption}{Assumption}

\newtheorem{remark}{Remark}

\newcommand{\bbR}{\mathbb{R}}

\newcommand{\ra}{\rightarrow}

\newcommand{\rlbrack}[1]{\left [ {#1} \right]}

\newcommand{\rlpar}[1]{\left ( {#1} \right)}

\title{Best Response Dynamics for Zero-Sum Dynamic Games \\ with Partial-Asymmetric Information}

\author{Yuxiang Guan \footnote{Ph.D. Candidate, Department of Mechanical Engineering; Email: yuxiang.guan@utdallas.edu (corresponding author).}}
\affil{The University of Texas at Dallas, Richardson, TX, 75080-3021}
\author{Iman Shames \footnote{Professor, Department of Electrical and Electronic Engineering; Email: iman.shames@unimelb.edu.au.}}
\affil{The University of Melbourne, Parkville, VIC 3010, Australia}
\author{Tyler H. Summers\footnote{Associate Professor, Department of Mechanical Engineering; Email: tyler.summers@utdallas.edu.}}
\affil{The University of Texas at Dallas, Richardson, TX, 75080-3021}

\begin{document}

\maketitle

\begin{abstract}
 This work studies a class of zero-sum stochastic linear quadratic dynamic games (LQDGs) under partial and asymmetric information. Information asymmetry introduces fundamental challenges related to \textit{belief representation} and \textit{theory of mind}, where players must impute belief states and estimates of other players to inform their strategies. Existing work highlights the difficulty of applying dynamic programming-like decomposition approach to these problems. An alternative approach based on \textit{best response dynamics} is proposed, which provides insights into belief representation and theory of mind challenges. Explicit expressions for each player's best response within the class of pure linear dynamic output feedback control strategies are derived, where the internal state dimension of each control is an integer multiple of the system state dimension. As players iteratively update their best responses, they form increasingly higher-order belief states, leading to infinite-dimensional internal states. However, numerical results reveal that the game's value converges after only a few iterations, suggesting that higher-order belief states provide vanishing benefit. This work further conducts numerical experiments to analyze the impact of asymmetric beliefs, belief orders, relative controllability and observability, and direct feed-through on a linear quadratic pursuit-evasion game's value.
\end{abstract}




\section{Introduction}
Pursuit-evasion games have been studied extensively in the literature, with critical applications in missile-target engagement \cite{shima2002time, shima2011optimal, weiss2016minimum}, air-combat tactics \cite{austin1990game, shinar1977analysis}, and space situational awareness \cite{pontani2009numerical, li2020saddle, vasal2021signaling}. These scenarios are characterized by adversarial interactions where the pursuer seeks to minimize capture distance or intercept time while the evader attempts to maximize these metrics or escape entirely. 
Dynamic game theory provides a framework for analyzing such interactions and designing equilibrium feedback strategies when each player has perfect information about the system state and model parameters \cite{bacsar1998dynamic,yuksel2024stochastic}. However, practical scenarios often involve partial and asymmetric information due to sensor limitations or differing information structures among players \cite{huang2021dynamic,shishika2021partial,cavalieri2014incomplete}. Under these conditions, traditional dynamic game theory is not directly applicable.

The fundamental challenges of \textit{belief representation} and \textit{theory of mind}, defined as the ability to attribute mental states (e.g., beliefs, intentions, or desires) to oneself and others \cite{scassellati2002theory}, arise in solving dynamic games with partial and asymmetric information. In these scenarios, players must impute the belief states and estimates of other players to inform their strategy. This leads to an infinite regress of higher-order beliefs, which renders even representing strategies and value functions difficult. A further challenge is that information asymmetry creates motivation for signaling, obfuscation, and deception, in which the players utilize their actions to affect the information of other players to indirectly optimize their objective. These challenges have been recognized in static game theory since the $1960$s, including in the seminal work of Harsanyi \cite{harsanyi1967games}. In the context of dynamic games, it means that control strategies may be infinite-dimensional, or equivalently, each agent must store and use its entire history of available information at each time to make optimal decisions. 


The difficulties of solving dynamic games with partial and asymmetric information were first identified and discussed by Isaacs in his classic book on differential games (Chapter 12)\cite{isaacsDifferential1999}. Since then, considerable interest has focused on differential games with linear dynamics and quadratic payoffs due to their analytical tractability and their relevance in the perturbation analysis of nonlinear pursuit-evasion games. In a stochastic dynamic game, a player can have one of three information structures: no information, imperfect information, or perfect information. The scenario where one player has perfect information and the other has imperfect information was explored in \cite{behn1968class,swarup2003linear}. Additionally, \cite{rhodes1969differential} examined cases where one of the two players has either no information or, under certain conditions, perfect information about the system state. These specific information structures allow the application of the \textit{separation theorem} from single-player optimal control to solve these games. However, this method cannot be directly applied to stochastic dynamic games where both players have distinct noisy measurements of the system state. This more general setting was studied in \cite{rhodes1969stochastic} for the zero-sum case and in \cite{chong1971stochastic} for the nonzero-sum case in which both players share the same objective function using constrained state estimators.
Formal solutions proposed in \cite{willman1969formal}, \cite{bagchi1981linear} address games where both players have imperfect information about the system state. These solutions assume that each player has perfect recall of their own information and the optimal strategies depend on all available information, meaning their control strategies may be infinite-dimensional.

More recent work \cite{nayyar2013common,gupta2014common,gupta2016dynamic,vasal2021signaling,sinha2016structured,ouyang2016dynamic,ouyang2025approach,tang2023dynamic,pachter2017lqg,altman2009stochastic,hambly2023linear} has sought to apply the sequential decomposition method to compute equilibrium strategies for stochastic dynamic games with partial and asymmetric information by defining specific common information structures. These equilibrium strategies and their associated beliefs can be categorized into two groups based on whether the common information-based beliefs are dependent on the players' strategies. The authors of \cite{nayyar2013common,gupta2014common,gupta2016dynamic} studied a class of dynamic games in which players' common information-based beliefs are independent of their strategies, meaning there is no signaling among them. This property enables them to use a dynamic programming-like backward recursive algorithm to find \textit{Markov perfect equilibria} of a transferred game, which are equivalently a class of Nash equilibria of the original game. A more general class of dynamic games where the players' common information-based beliefs are strategy-dependent, that is, signaling is present, was explored in \cite{ouyang2016dynamic,tang2023dynamic,vasal2021signaling,sinha2016structured}. A sequential decomposition method was used to compute \textit{perfect Bayesian equilibria} for these games. All the work described above assumes the common information structure of the game and seeks a sequential decomposition method to compute equilibrium strategies. However, for general games lacking this common information structure, to the best of our knowledge, no known sequential decomposition method exists. Consequently, we do not provide a direct numerical comparison with common information-based sequential decomposition approaches in this paper, since these methods are not well defined for the private, noisy measurement setting studied here, and adapting them to this asymmetric structure remains an open problem in its own right. Our best response dynamics, by contrast, only require solving standard LQG subproblems at each iteration, so they apply directly without assuming any common information structure.

Recently, in \cite{guan2025best}, we analyzed the performance of best response dynamics for solving a two-player zero-sum linear quadratic Gaussian (LQG) dynamic game with partial and asymmetric information (where the state dynamics are linear, the cost function is quadratic, and each player receives only a noisy linear measurement of the state at each time) without the common information assumption. The best response dynamics work by initializing player $2$'s strategy to a zero strategy; then player $1$'s best response for this game is an LQG controller, combining a linear quadratic regulator (LQR) with a Kalman filter to estimate the underlying system state from their output measurements. This controller has an internal dimension equal to the underlying system state dimension, and the conditional state distribution represents a ($0$-th order) belief state. Then, fixing player $1$'s strategy, player $2$'s best response is an augmented LQG controller in which the Kalman filter estimates not only the underlying system state but also the state estimate of player $1$ (a $1$st-order belief state). With each successive best response between these two players, they form increasingly higher-order belief states, leading to an infinite regress where their internal state dimension increases towards infinity. However, we observed in extensive numerical experiments that the game's value converges after just a few iterations due to the significant difficulty in controlling and observing both players' higher-order belief dynamics.

In this work, we extend the analysis presented in \cite{guan2025best} along two directions: the analytical derivation of the best response and the scope of the numerical study. Analytically, we incorporate direct feed-through of opposing players' inputs into each player's measurement model, which is much more challenging than in single-player settings because the feed-through effect cannot be simply subtracted out due to the information asymmetry, and each player has an incentive to use their actions to affect the opposing player's measurements. This is not a direct extension of the solution in \cite{guan2025best}; each player's augmented filter, belief-state recursion, and resulting controller must be rederived to account for the opponent's ability to corrupt measurements. Numerically, we move beyond the low-dimensional, randomly generated open-loop stable systems used in \cite{guan2025best} to a physically motivated linear quadratic pursuit-evasion game, instantiated with a spacecraft relative-motion model, and we generalize the Gramian and Hankel singular value analysis of higher-order belief dynamics to a broader set of randomly generated higher-dimensional systems. Our main contributions are threefold:
\begin{enumerate}
    \item We formulate a class of two-player zero-sum stochastic LQDGs where each player's measurement contains the opponent's input through an unobservable direct feed-through term that cannot be subtracted, instantiated as a spacecraft pursuit-evasion game where each player operates two independent subsystems, a maneuvering input affecting the relative kinematics and an electronic countermeasure (ECM) input corrupting only the opponent's sensor measurement.
    \item We derive explicit alternating best response expressions within the class of pure linear dynamic output feedback strategies, showing that information asymmetry and the unobservable feed-through term force each player's internal belief state dimension to grow by an integer multiple of the system state dimension at each iteration.
    \item We investigate numerically the impacts of asymmetric beliefs, belief orders, relative controllability and observability, and ECM effectiveness on the pursuit-evasion game's value and chase trajectories, with Gramian eigenvalues and Hankel singular values decay in higher-order belief dynamics providing the mechanistic explanation for observed convergence within a few iterations.
\end{enumerate}

The remainder of this paper is structured as follows. In Section \ref{sec:problem_formulation}, we formulate a class of infinite horizon two-player zero-sum stochastic LQDGs with partial and asymmetric information with an instance of linear quadratic pursuit-evasion games. Section \ref{sec:best_response_adver_inputs} presents explicit expressions for each player's best response, taking into account direct feed-through in their measurements. In Section \ref{sec:experimental_results}, we analyze the effects of players' asymmetric beliefs, belief orders, relative controllability, observability, and measurement manipulations on the game's value. We conclude our results and discussions in Section \ref{sec:conclusion}.

\section{Problem Formulation}\label{sec:problem_formulation}
\subsection{Zero-Sum Games with Direct Feed-Through}\label{sec:general-zeo-sum-games}
We consider a class of two-player zero-sum LQDGs satisfying the following dynamics
\begin{equation}\label{eqn:zero_sum_dynamics}
    \mathbf{x}_{t+1} = \mathbf{A} \mathbf{x}_t + \mathbf{B}^{(1)} \mathbf{u}_t^{(1)} + \mathbf{B}^{(2)} \mathbf{u}_t^{(2)} + \mathbf{w}_t,
\end{equation}
where $\mathbf{x}_t \in \bbR^n$ represents the system state at time $t$ with an initial random value $\mathbf{x}_0 \sim \mathcal{N}(\mathbf{\bar x}_0, \mathbf{X}_0)$. The control actions of player~1 and player~2 at time $t$ are $\mathbf{u}_t^{(1)} \in \bbR^{m_1}$ and $\mathbf{u}_t^{(2)} \in \bbR^{m_2}$. The process noise $\mathbf{w}_t \in \bbR^n$ is an independent and identically distributed random vector with distribution $\mathcal{N}(0, \mathbf{W})$. Furthermore, the system and input matrices are $\mathbf{A} \in \bbR^{n \times n}$, $\mathbf{B}^{(1)} \in \bbR^{n \times m_1}$, and $\mathbf{B}^{(2)} \in \bbR^{n \times m_2}$. 

In this work, we explore a setting characterized by partial and asymmetric information, where players lack precise knowledge of the system states and the opposing player's actions. This uncertainty creates an incentive for each player to disrupt the opponent’s measurements and thereby influence their beliefs about the true system states. To capture this dynamic, we incorporate the direct feed-through of the opponent's inputs into each player's measurement model, reflecting both the intent and the capability to interfere with the opponent's measurements through the following model:
\begin{align}\label{eqn:zero-sum-measurement}
    \begin{split}
        \mathbf{y}_t^{(1)} &= \mathbf{C}^{(1)} \mathbf{x}_t + \mathbf{D}^{(2)} \mathbf{u}_t^{(2)} + \mathbf{v}_t^{(1)}, \\
        \mathbf{y}_t^{(2)} &= \mathbf{C}^{(2)} \mathbf{x}_t + \mathbf{D}^{(1)} \mathbf{u}_t^{(1)} + \mathbf{v}_t^{(2)},
    \end{split}
\end{align}
where $\mathbf{y}_t^{(1)} \in \bbR^{p_1}$ and $\mathbf{y}_t^{(2)} \in \bbR^{p_2}$ are the observed outputs of player~1 and player~2, $\mathbf{v}_t^{(1)} \in \bbR^{p_1}$ and $\mathbf{v}_t^{(2)} \in \bbR^{p_2}$ are independent and identically distributed random vectors with distributions $\mathcal{N}(0, \mathbf{V}^{(1)})$ and $\mathcal{N}(0, \mathbf{V}^{(2)})$. The output matrices are $\mathbf{C}^{(1)} \in \bbR^{p_1 \times n}$ and $\mathbf{C}^{(2)} \in \bbR^{p_2 \times n}$. The feed-through matrices are $\mathbf{D}^{(1)} \in \bbR^{p_2 \times m_1}$ and $\mathbf{D}^{(2)} \in \bbR^{p_1 \times m_2}$.

To characterize the objective of each player in classic zero-sum games, we introduce a quadratic functional of $\mathbf{x}_t$, $\mathbf{u}_t^{(1)}$, and $\mathbf{u}_t^{(2)}$ over an infinite horizon as
\begin{equation}\label{eqn:obj_fun}
    J = \lim_{T \ra \infty} \frac{1}{T}\mathbf{E}_{\mathbf{x}_0, \mathbf{w}_t, \mathbf{v}_t^{(1)}, \mathbf{v}_t^{(2)}} \Biggl[\sum_{t=0}^{T-1} \mathbf{x}_t^\top \mathbf{Q} \mathbf{x}_t + \rlpar{\mathbf{u}_t^{(1)}}^\top \mathbf{R}^{(1)} \mathbf{u}_t^{(1)} + \rlpar{\mathbf{u}_t^{(2)}}^\top \mathbf{R}^{(2)} \mathbf{u}_t^{(2)}\Biggr],
\end{equation}
where $\mathbf{Q} \in \bbR^{n \times n} \succeq 0$ is the state penalty matrix, $\mathbf{R}^{(1)} \in \bbR^{m_1 \times m_1} \succ 0$, and $\mathbf{R}^{(2)} \in \bbR^{m_2 \times m_2} \prec 0$ are the penalty matrices for player~1 and player~2's control actions. To ensure the game's upper value is bounded, we assume the penalty $\mathbf{R}^{(2)}$ is sufficiently large. In classic zero-sum games, player~1 seeks strategies to minimize the objective function \eqref{eqn:obj_fun} while player~2 seeks strategies to maximize it. 

Since the players do not know the exact value of the states, they must in general make decisions based on a history of all available information. Accordingly, we define the information states $\mathbf{I}_t^{(i)}$ ($i = 1, 2$) for $t=0,1,\ldots$ as
\begin{align*}
    \begin{split}
        \mathbf{I}_t^{(1)} &= \left( \mathbf{y}_0^{(1)}, \mathbf{y}_1^{(1)}, \ldots, \mathbf{y}_t^{(1)}, \mathbf{u}_0^{(1)}, \mathbf{u}_1^{(1)}, \ldots, \mathbf{u}_{t-1}^{(1)} \right), \\
        \mathbf{I}_t^{(2)} &= \left( \mathbf{y}_0^{(2)}, \mathbf{y}_1^{(2)}, \ldots, \mathbf{y}_t^{(2)}, \mathbf{u}_0^{(2)}, \mathbf{u}_1^{(2)}, \ldots, \mathbf{u}_{t-1}^{(2)} \right).
    \end{split}
\end{align*}
The information setting defined above means that each player observes only their own noisy measurement $\mathbf{y}_t^{(i)}$ and knows their own past controls $\mathbf{u}_t^{(i)}$, but has no access to the opponent's measurements and controls. Moreover, each player's measurement $\mathbf{y}_t^{(i)}$ is distinct, arising from distinct measurement noises $\mathbf{v}_t^{(1)}\neq \mathbf{v}_t^{(2)}$, distinct feed-through terms $\mathbf{D}^{(1)}\mathbf{u}_t^{(1)} \neq \mathbf{D}^{(2)}\mathbf{u}_t^{(2)}$, and/or distinct output matrices $\mathbf{C}^{(1)} \neq \mathbf{C}^{(2)}$.
The players seek strategies designed by information feedback functions that optimize the objective function \eqref{eqn:obj_fun}. More specifically, we consider a class of pure linear dynamic output feedback strategies where internal state dimension of each control strategy is an integer multiple of the system state dimension, i.e., we seek $\mathbf{\pi}^{(1)}$ with $\mathbf{u}_t^{(1)} = \mathbf{\pi}^{(1)}(\mathbf{z}_t^{(1)})$ and $\mathbf{\pi}^{(2)}$ with $\mathbf{u}_t^{(2)} = \mathbf{\pi}^{(2)}(\mathbf{z}_t^{(2)})$ where $\mathbf{z}_t^{(1)} \in \bbR^{n_k^{(1)} n}$ and $\mathbf{z}_t^{(2)} \in \bbR^{n_k^{(2)} n}$ are state estimates (or belief states) with  $n_k^{(1)}, n_k^{(2)} \in \{1, 2, \ldots\}$ ($k$ is the number of best response iterations) produced by filters of the form
\begin{align}\label{eqn:belief-dynamics}
    \begin{split}
        \mathbf{z}_{t+1}^{(1)} &= \mathbf{A}^{(1)} \mathbf{z}_t^{(1)} + \bar{\mathbf{B}}^{(1)} \mathbf{u}_t^{(1)} + \mathbf{L}^{(1)} \left( \mathbf{y}_t^{(1)} - \bar{\mathbf{C}}^{(1)} \mathbf{z}_t^{(1)} \right), \\
        \mathbf{z}_{t+1}^{(2)} &= \mathbf{A}^{(2)} \mathbf{z}_t^{(2)} + \bar{\mathbf{B}}^{(2)} \mathbf{u}_t^{(2)} + \mathbf{L}^{(2)} \left( \mathbf{y}_t^{(2)} - \bar{\mathbf{C}}^{(2)} \mathbf{z}_t^{(2)} \right).
    \end{split}
\end{align}
where $\mathbf{A}^{(i)} \in \bbR^{n_k^{(i)} n \times n_k^{(i)} n}$, $\bar{\mathbf{B}}^{(i)} \in \bbR^{n_k^{(i)} n \times m_i}$, and $\bar{\mathbf{C}}^{(i)} \in \bbR^{p_i \times n_k^{(i)} n}$ are the system, input, and output matrices of the players' higher-order belief dynamics, $\mathbf{L}^{(i)} \in \bbR^{n_k^{(i)} n \times p_i}$ denotes the higher-order state estimator. The dimension of $\mathbf{z}_t^{(i)}$ increases with each iteration because player $i$'s belief state must account for the opposing player's belief state, which itself grows by $n$ at each iteration.

\begin{assumption}\label{assmp1}
We assume that all system parameters (comprising the system and input matrices $\mathbf{A}$, $\mathbf{B}^{(1)}$, $\mathbf{B}^{(2)}$, the output matrices $\mathbf{C}^{(1)}$ and $\mathbf{C}^{(2)}$, the feed-through matrices $\mathbf{D}^{(1)}$ and $\mathbf{D}^{(2)}$, the process noise covariance $\mathbf{W}$, the measurement noise covariances $\mathbf{V}^{(1)}$ and $\mathbf{V}^{(2)}$, and the initial state mean $\mathbf{\bar{x}}_0$ and covariance $\mathbf{X}_0$) and all penalty parameters (comprising the state penalty matrix $\mathbf{Q}$, player~1's penalty matrix $\mathbf{R}^{(1)}$, and player~2's penalty matrix $\mathbf{R}^{(2)}$) are common knowledge across all players.
\end{assumption}

\begin{remark}
    Assumption \ref{assmp1} is standard in the literature on dynamic games \cite{bacsar1998dynamic}, which enables each player to independently reason about their best responses and those of the other player to determine optimal strategies. We acknowledge that this assumption has limits in some practical scenarios, where players may be uncertain about each other's system or penalty parameters; in this case, the game becomes one of incomplete information in the sense of \cite{harsanyi1967games}, requiring each player to form beliefs over the opponent's unknown parameters and optimally respond against those beliefs. One natural way to incorporate such uncertainty is to replace an unknown parameter, for example the feed-through matrix $\mathbf{D}^{(2)}$ representing player~1's uncertainty about whether player~2 can affect player~1's measurement, with a finite set of hypothesized values, each assigned a prior probability known to both players. Player~1 would then maintain a separate belief state for each hypothesized value, update a posterior over the hypotheses from its measurements, and respond with a policy that mixes the per-hypothesis best responses according to this posterior.

    This extension is significantly harder than the complete-information setting studied in this paper, since each player would need to track a separate belief hierarchy for every hypothesized parameter value in addition to the evolving posterior over those hypotheses, and no general closed-form solution to this class of problems exists even in the single-player LQG setting. We regard this as a separate and substantial research problem in its own right, and leave it to future work.
\end{remark}

 In this work, we are interested in finding a saddle point (Nash) equilibrium policy pair ($\pi^{(1)*}$, $\pi^{(2)*}$) where
\begin{equation*}
    J(\pi^{(1)*}, \pi^{(2)}) \leq J(\pi^{(1)*}, \pi^{(2)*}) \leq J(\pi^{(1)}, \pi^{(2)*}), \quad \forall \pi^{(1)}, \pi^{(2)}.
\end{equation*}
This means that no player can unilaterally improve their cost by deviating from their equilibrium policy.

\subsection{Two-Player Pursuit-Evasion Games}\label{sec:two_player_pe_game}
A two-player linear quadratic pursuit-evasion game is a natural instance of the zero-sum stochastic LQDG formulation in Section \ref{sec:general-zeo-sum-games}, in which the pursuer takes the role of the minimizing player and the evader takes the role of the maximizing player. This pursuit-evasion framework accommodates a broad class of aerospace and aeronautical engagement scenarios: the underlying dynamics can represent linearized missile-target terminal homing kinematics \cite{shima2002time,shima2011optimal}, close-range air combat maneuvering between fighter aircraft (dogfight), or spacecraft proximity operations in orbit, depending on the physical setting of interest. In this work, we consider a spacecraft pursuit-evasion scenario governed by orbital relative motion dynamics, which provides a physically grounded and analytically tractable instantiation of the general framework. 

\subsubsection{Dynamical Model}
We consider a pursuit-evasion engagement between two spacecraft operating in proximity in low Earth orbit, representative of space situational awareness and orbital inspection scenarios \cite{pontani2009numerical, li2020saddle}. To formulate the relative dynamics, we use the local-vertical local-horizontal (LVLH) frame centered on the evader spacecraft, with its origin fixed to the evader. In this frame, the $x$-axis coincides with the orbital radius vector, positive outward; the $z$-axis is oriented toward the angular momentum vector; and the $y$-axis completes the right-handed Cartesian frame.

Under the assumption of a near-circular reference orbit and small deviations from a nominal condition, the linearized relative equations of motion in the LVLH frame are given by the Clohessy-Wiltshire (CW) equations:
\begin{align*}
    \ddot{x} - 3wx - 2w\dot{y} &= u_x^{(p)} - u_x^{(e)} + d_x^{(p)} - d_x^{(e)}, \\
    \ddot{y} + 2w\dot{x} &= u_y^{(p)} - u_y^{(e)} + d_y^{(p)} - d_y^{(e)}, \\
    \ddot{z} + w^2 z &= u_z^{(p)} - u_z^{(e)} + d_z^{(p)} - d_z^{(e)},
\end{align*}
where $w = \sqrt{\mu/a^3} \ \mathrm{rad/s}$ is the orbital rate of the target body, $a$ is the radius of the target body's circular orbit, $\mu$ is the standard gravitational parameter; $(x,\ y,\ z) \ (\mathrm{m})$ denote the relative position of the pursuer with respect to the evader in the radial, along-track, and cross-track directions; $\mathbf{u}^{(p)}=[u_x^{(p)},\ u_y^{(p)},\ u_z^{(p)}]^\top\in\mathbb{R}^3$ and $\mathbf{u}^{(e)} = [u_x^{(e)},\ u_y^{(e)},\ u_z^{(e)}]^\top\in\mathbb{R}^3$ are the thrust acceleration commands $(\mathrm{m/s^2})$ of the pursuer and evader; ${\mathbf{d}^{(p)}}=[d_x^{(p)},\ d_y^{(p)},\ d_z^{(p)}]^\top$ and ${\mathbf{d}^{(e)}}=[d_x^{(e)},\ d_y^{(e)},\ d_z^{(e)}]^\top$ are zero-mean Gaussian disturbances with covariance $\mathbf{W}_d^{(p)}$ and $\mathbf{W}_d^{(e)}$ representing unmodeled perturbations such as atmospheric drag and solar radiation pressure.

Defining the relative state vector $\mathbf{x}=[x,\ y,\ z,\ \dot{x},\ \dot{y},\ \dot{z}]^\top\in\mathbb{R}^6$, where the first three components are relative positions $(\mathrm{m})$ and the last three are relative velocities $(\mathrm{m/s})$, the CW equations take the continuous-time state-space form:
\begin{equation}
    \dot{\mathbf{x}} = \mathbf{A}_c \mathbf{x} + \mathbf{B}_c^{(p)} \mathbf{u}^{(p)} + \mathbf{B}_c^{(e)} \mathbf{u}^{(e)} + \mathbf{G}_c^{(p)} \mathbf{d}^{(p)} + \mathbf{G}_c^{(e)} \mathbf{d}^{(e)},
\end{equation}
where 
\begin{equation*}
    \mathbf{A}_c = \begin{bmatrix}
        0 & 0 & 0 & 1 & 0 & 0 \\
        0 & 0 & 0 & 0 & 1 & 0 \\
        0 & 0 & 0 & 0 & 0 & 1 \\
        3w^2 & 0 & 0 & 0 & 2w & 0 \\
        0 & 0 & 0 & -2w & 0 & 0 \\
        0 & 0 & -w^2 & 0 & 0 & 0
    \end{bmatrix}, \quad \mathbf{B}_c^{(p)} = \mathbf{G}_c^{(p)} = \begin{bmatrix}
        0 & 0 & 0 \\
        0 & 0 & 0 \\
        0 & 0 & 0 \\
        1 & 0 & 0 \\
        0 & 1 & 0 \\
        0 & 0 & 1
    \end{bmatrix}, \quad \mathbf{B}_c^{(e)} = -\mathbf{B}_c^{(p)}, \quad \mathbf{G}_c^{(e)} = -\mathbf{G}_c^{(p)}.
\end{equation*}

To discretize the continuous-time CW equations using the Zero-Order Hold (ZOH) method with sampling period $\Delta t$, defining the parameter $\phi = w\Delta t$. The discrete-time state-space model is
\begin{equation}
    \mathbf{x}_{t+1} = \mathbf{A}\mathbf{x}_t + \mathbf{B}^{(p)}\mathbf{u}_t^{(p)} + \mathbf{B}^{(e)}\mathbf{u}_t^{(e)} + \tilde{\mathbf{w}}_t,
\end{equation}
where
\begin{equation*}
    \mathbf{A} = \begin{bmatrix}
        4-3\cos{\phi} & 0 & 0 & \frac{\sin{\phi}}{w} & \frac{2}{w}(1-\cos{\phi}) & 0 \\
        -6(\phi - \sin{\phi}) & 1 & 0 & -\frac{2}{w}(1-\cos{\phi}) & \frac{4\sin{\phi}-3\phi}{w} & 0 \\
        0 & 0 & \cos{\phi} & 0 & 0 & \frac{\sin{\phi}}{w} \\
        3w\sin{\phi}  & 0 & 0 & \cos{\phi} & 2\sin{\phi} & 0 \\
        -6w(1-\cos{\phi}) & 0 & 0 & -2\sin{\phi} & 4\cos{\phi}-3 & 0 \\
        0 & 0 & -w\sin{\phi} & 0 & 0 & \cos{\phi}
    \end{bmatrix},
\end{equation*}
\begin{equation*}
    \mathbf{B}^{(p)} = \mathbf{G}^{(p)} = \begin{bmatrix}
    \frac{1-\cos{\phi}}{w^2} & \frac{2(\phi-\sin{\phi})}{w^2} & 0 \\
    -\frac{2(\phi-\sin{\phi})}{w^2} & \frac{8(1-\cos{\phi})-3\phi^2}{2w^2} & 0 \\
    0 & 0 & \frac{1-\cos{\phi}}{w^2} \\
    \frac{\sin{\phi}}{w} & \frac{2(1-\cos{\phi})}{w} & 0 \\
    -\frac{2(1-\cos{\phi})}{w} & \frac{4\sin{\phi}-3\phi}{w} & 0 \\
    0 & 0 & \frac{\sin{\phi}}{w}
\end{bmatrix}, \quad \mathbf{B}^{(e)} = -\mathbf{B}^{(p)}, \quad \mathbf{G}^{(e)} = -\mathbf{G}^{(p)},
\end{equation*}
and $\tilde{\mathbf{w}}_t=\mathbf{G}^{(p)}\mathbf{d}_t^{(p)}+\mathbf{G}^{(e)}\mathbf{d}_t^{(e)} \sim \mathcal{N}(\mathbf{0}, \tilde{\mathbf{W}})$ is the discretized process noise with $\tilde{\mathbf{W}} = \mathbf{G}^{(p)} \mathbf{W}_d^{(p)} \mathbf{G}^{{(p)}\top} + \mathbf{G}^{(e)} \mathbf{W}_d^{(e)} \mathbf{G}^{{(e)}\top}$.

\subsubsection{Measurement Model}
Since neither player directly observes the full relative state, each player receives noisy measurements of the relative position through an onboard ranging sensor. In addition to maneuvering via thrust input $\mathbf{u}_t^{(i)}$ ($i=p,e$), each player may operate an ECM subsystem that deliberately corrupts the opponent's sensor measurement. Denoting the pursuer's and evader's ECM signals by $\mathbf{s}_t^{(p)} \in\mathbb{R}^{q_p}$ and $\mathbf{s}_t^{(e)} \in \mathbb{R}^{q_e}$ respectively. The measurement received by each player is
\begin{align}\label{eq:sensor_ecm}
    \begin{split}
        \mathbf{y}_t^{(p)} &= \mathbf{C}^{(p)} \mathbf{x}_t + \mathbf{D}^{(e)} \mathbf{s}_t^{(e)} + \mathbf{v}_t^{(p)}, \\
        \mathbf{y}_t^{(e)} &= \mathbf{C}^{(e)} \mathbf{x}_t + \mathbf{D}^{(p)} \mathbf{s}_t^{(p)} + \mathbf{v}_t^{(e)},
    \end{split}
\end{align}
where $\mathbf{C}^{(p)} = \mathbf{C}^{(e)} = [\mathbf{I}_3,\ \mathbf{0}_3] \in \mathbb{R}^{3 \times 6}$ selects only the relative position states, reflecting a realistic sensor model in which relative position is available but relative velocity must be estimated from the filter. The matrices $\mathbf{D}^{(p)} \in\mathbb{R}^{3 \times q_p}$ and $\mathbf{D}^{(e)} \in\mathbb{R}^{3 \times q_e}$ are the jamming effectiveness matrices quantifying how strong each player's ECM signal corrupts the opponent's measurement -- representative of radar jamming, GPS spoofing, optical decoys, or laser dazzling of the opponent's sensor. The measurement noises $\mathbf{v}_t^{(p)} \sim \mathcal{N}(\mathbf{0}, \mathbf{V}^{(p)})$ and $\mathbf{v}_t^{(e)} \sim \mathcal{N}(\mathbf{0}, \mathbf{V}^{(e)})$ are independent across players over time.

\subsubsection{Connection to Zero-Sum Games}
In a classic pursuit-evasion game, the pursuer aims to minimize the relative state $\mathbf{x}_t$ to capture the evader, while the evader seeks to maximize it to avoid capture. If capture is inevitable, the evader aims to maximize the time until capture. Since each player operates two independent subsystems, maneuvering via thrust input $\mathbf{u}_t^{(i)}$ and sensor interference via ECM input $\mathbf{s}_t^{(i)}$, we define the augmented control inputs
\begin{equation*}
    \tilde{\mathbf{u}}_t^{(p)} = \begin{bmatrix}
        \mathbf{u}_t^{(p)} \\
        \mathbf{s}_t^{(p)}
    \end{bmatrix} \in \mathbb{R}^{m_p + q_p}, \quad\quad \tilde{\mathbf{u}}_t^{(e)} = \begin{bmatrix}
        \mathbf{u}_t^{(e)} \\
        \mathbf{s}_t^{(e)}
    \end{bmatrix} \in \mathbb{R}^{m_e + q_e}.
\end{equation*}

The framework of Section \ref{sec:general-zeo-sum-games} can be applied to compute Nash equilibrium for this pursuit-evasion game by identifying:
\begin{equation*}
    \mathbf{B}^{(1)} = \tilde{\mathbf{B}}^{(p)}, \quad \mathbf{B}^{(2)} = \tilde{\mathbf{B}}^{(e)}, \quad \mathbf{C}^{(1)} = \mathbf{C}^{(p)}, \quad \mathbf{C}^{(2)} = \mathbf{C}^{(e)}, \quad \mathbf{D}^{(1)} = \tilde{\mathbf{D}}^{(p)}, \quad \mathbf{D}^{(2)} = \tilde{\mathbf{D}}^{(e)},
\end{equation*}
\begin{equation*}
    \mathbf{w}_t = \tilde{\mathbf{w}}_t \sim \mathcal{N}(\mathbf{0}, \tilde{\mathbf{W}}), \quad \mathbf{v}_t^{(1)} = \mathbf{v}_t^{(p)}, \quad \mathbf{v}_t^{(2)} = \mathbf{v}_t^{(e)},
\end{equation*}
where the augmented system matrices are
\begin{equation*}
    \tilde{\mathbf{B}}^{(p)} = \begin{bmatrix}
        \mathbf{B}^{(p)} & \mathbf{0}
    \end{bmatrix}, \quad \tilde{\mathbf{B}}^{(e)} = \begin{bmatrix}
        \mathbf{B}^{(e)} & \mathbf{0}
    \end{bmatrix}, \quad \tilde{\mathbf{D}}^{(p)} = \begin{bmatrix}
        \mathbf{0} & \mathbf{D}^{(p)}
    \end{bmatrix}, \quad \tilde{\mathbf{D}}^{(e)} = \begin{bmatrix}
        \mathbf{0} & \mathbf{D}^{(e)}
    \end{bmatrix}.
\end{equation*}
The penalty matrices in Section \ref{sec:general-zeo-sum-games} are correspondingly identified as
\begin{equation*}
    \mathbf{R}^{(1)} = \begin{bmatrix}
        \mathbf{R}_u^{(p)} & \mathbf{0} \\
        \mathbf{0} & \mathbf{R}_s^{(p)}
    \end{bmatrix} \succ 0, \quad\quad \mathbf{R}^{(2)} = \begin{bmatrix}
        \mathbf{R}_u^{(e)} & 0 \\
        \mathbf{0} & \mathbf{R}_s^{(e)}
    \end{bmatrix} \prec 0,
\end{equation*}
where $\mathbf{R}_u^{(p)}$, $\mathbf{R}_s^{(p)} \succ 0$ penalize the pursuer's maneuvering effort and ECM power independently, and $\mathbf{R}_u^{(e)}$, $\mathbf{R}_s^{(e)} \prec 0$ are the corresponding matrices for the evader. The analytical insights gained from this linear quadratic setting provide a foundation for extending the proposed best response dynamics approach to nonlinear games through successive linearization implementations.


\section{Best Response Dynamics with Direct feed-through in Players' Measurements}\label{sec:best_response_adver_inputs}
The key idea of using best response dynamics to solve the zero-sum stochastic LQDGs described above is that, at each iteration, each player solves a standard LQG problem for an augmented system. Specifically, player $2$'s strategy is initialized to zero, so that player $1$ first solves a standard LQG problem. Player $2$ then treats player $1$'s resulting strategy as fixed and solves its own LQG problem, now over an augmented state that incorporates player $1$'s zeroth-order belief state (estimate about the system state). In subsequent iterations, each player must form beliefs not only about the system state but also about the opponent's belief state from the previous iteration. This leads to an augmented LQG problem for each player to solve at each iteration.

Although each player individually solves an LQG subproblem at each iteration, two features fundamentally distinguish this framework from the standard single-player LQG problem. First, each player's measurement \eqref{eqn:zero-sum-measurement} contains the opponent's input through the direct feed-through term, for example, $\mathbf{D}^{(1)}\mathbf{u}_t^{(1)}$ enters player $2$'s measurement, which player $2$ cannot observe and therefore cannot subtract as in single-player case. This term must instead be estimated, coupling estimation and control across players in a way that has no single-player analogue. Second, each player's filter dimension grows by $n$ at every iteration, rather than remaining fixed at the system state dimension $n$ as in single-player LQG, leading to infinite-dimensional strategies in the limit. In practice, however, Section \ref{sec:experimental_results} shows that convergence occurs after only a few iterations.

\subsection{LQG Control of Player 1 (pursuer/minimizer)}\label{sec:lqg_pursuer}
We first initialize a zero strategy for player $2$ (where $\mathbf{u}_t^{(2)} = 0 \ \forall t$). This initialization choice is not strictly required when the dynamic system is open-loop stable: one may instead initialize with a zero strategy for player $1$, or with any pair of random stabilizable strategies for both players, and the best response dynamics converge to the same game value regardless. For open-loop unstable or marginally stable systems, however, initializing with player $2$'s zero strategy is preferred, since it reduces the game to a 
standard LQG problem for player $1$ and guarantees that the upper value of the game remains bounded from the first iteration.

Then player $1$ solves a standard LQG control problem with dynamics and measurement
\begin{align}\label{eqn:sys_player1}
\begin{split}
        \mathbf{x}_{t+1} &= \mathbf{A} \mathbf{x}_t + \mathbf{B}^{(1)} \mathbf{u}_t^{(1)} + \mathbf{w}_t, \\
        \mathbf{y}_t^{(1)} &= \mathbf{C}^{(1)} \mathbf{x}_t + \mathbf{v}_t^{(1)},
\end{split}
\end{align}
and objective 
\begin{equation*}
    J^{(1)}  = \lim_{T \ra \infty} \frac{1}{T} \mathbf{E}_{\mathbf{x}_0, \mathbf{w}_t, \mathbf{v}_t^{(1)}} \rlbrack{\sum_{t=0}^{T-1} \mathbf{x}_t^\top \mathbf{Q} \mathbf{x}_t + \rlpar{\mathbf{u}_t^{(1)}}^\top \mathbf{R}^{(1)} \mathbf{u}_t^{(1)}}.
\end{equation*}

Player $1$ employs a steady-state Kalman filter to form an estimate $\mathbf{z}_t^{(1)} \in \mathbb{R}^n$ of the state given by
\begin{align}
    \mathbf{z}_{t+1}^{(1)} = \mathbf{A} \mathbf{z}_t^{(1)} + \mathbf{B}^{(1)} \mathbf{u}_t^{(1)} + \mathbf{L}^{(1)} \rlpar{\mathbf{y}_t^{(1)} - \mathbf{C}^{(1)} \mathbf{z}_t^{(1)}},
\end{align}
with 
\begin{equation}
    \mathbf{L}^{(1)} = \mathbf{A} \mathbf{\Sigma}^{(1)} \rlpar{\mathbf{C}^{(1)}}^\top \rlpar{\mathbf{V}^{(1)} + \mathbf{C}^{(1)} \mathbf{\Sigma}^{(1)} \rlpar{\mathbf{C}^{(1)}}^\top}^{-1},
\end{equation}
where the optimal estimator gain $\mathbf{L}^{(1)}$ utilizes the steady-state estimation error covariance matrix $\mathbf{\Sigma}^{(1)}$ that solves the algebraic Riccati equation
\begin{align}\label{eqn:Sigma1}
        \mathbf{\Sigma}^{(1)} = \mathbf{W} + \mathbf{A} \mathbf{\Sigma}^{(1)} \mathbf{A}^\top - \mathbf{A} \mathbf{\Sigma}^{(1)} \rlpar{\mathbf{C}^{(1)}}^\top
         \rlpar{\mathbf{V}^{(1)} + \mathbf{C}^{(1)} \mathbf{\Sigma}^{(1)} \rlpar{\mathbf{C}^{(1)}}^\top}^{-1} \mathbf{C}^{(1)} \mathbf{\Sigma}^{(1)} \mathbf{A}^\top.
\end{align}
The steady-state estimation error covariance is given by $\mathbf{\Sigma}^{(1)} = \lim_{t\rightarrow \infty} \mathbf{\Sigma}_t^{(1)}$, where $\mathbf{\Sigma}_t^{(1)} = \mathbf{E} [ \mathbf{e}_t^{(1)} (\mathbf{e}_t^{(1)})^\top]$ and the estimation error $\mathbf{e}_t^{(1)} \coloneqq \mathbf{x}_t - \mathbf{z}_t^{(1)} $ has dynamics
\begin{align*}
    \begin{split}
        \mathbf{e}_{t+1}^{(1)} &= \mathbf{x}_{t+1} - \mathbf{z}_{t+1}^{(1)} \\
                  &= \rlpar{\mathbf{A} - \mathbf{L}^{(1)} \mathbf{C}^{(1)}} \mathbf{e}_t^{(1)} + \mathbf{w}_t - \mathbf{L}^{(1)} \mathbf{v}_t^{(1)}.
    \end{split}
\end{align*}

The optimal strategy of player $1$ is the linear state estimate feedback 
\begin{align}\label{eqn:u_t_1}
    \mathbf{u}_t^{(1)} = \mathbf{K}^{(1)} \mathbf{z}_t^{(1)},
\end{align}
where the gain matrix $\mathbf{K}^{(1)} \in \mathbb{R}^{m_1 \times n}$ is given by
\begin{align}\label{eqn:K1}
    \mathbf{K}^{(1)} = -\rlpar{\mathbf{R}^{(1)} + \rlpar{\mathbf{B}^{(1)}}^\top \mathbf{P}^{(1)} \mathbf{B}^{(1)}}^{-1} \rlpar{\mathbf{B}^{(1)}}^\top \mathbf{P}^{(1)} \mathbf{A},
\end{align}
where the cost matrix $\mathbf{P}^{(1)}$ solves the algebraic Riccati equation
\begin{align}\label{eqn:P1}
    \mathbf{P}^{(1)} = \mathbf{Q} + \mathbf{A}^\top \mathbf{P}^{(1)} \mathbf{A} - \mathbf{A}^\top \mathbf{P}^{(1)} \mathbf{B}^{(1)} \rlpar{\mathbf{R}^{(1)} + \rlpar{\mathbf{B}^{(1)}}^\top \mathbf{P}^{(1)} \mathbf{B}^{(1)}}^{-1} \rlpar{\mathbf{B}^{(1)}}^\top \mathbf{P}^{(1)} \mathbf{A}.
\end{align}
There exist unique solutions $\mathbf{P}^{(1)}$ and $\mathbf{\Sigma}^{(1)}$ of the Riccati equations if the system \eqref{eqn:sys_player1} is stabilizable and detectable. 

The closed-loop state and state estimate dynamics are
\begin{align*} \nonumber
     \begin{bmatrix}
         \mathbf{x}_{t+1} \\ \mathbf{z}_{t+1}^{(1)}
     \end{bmatrix} &= \underbrace{\begin{bmatrix}
        \mathbf{A} & \mathbf{B}^{(1)} \mathbf{K}^{(1)} \\
        \mathbf{L}^{(1)} \mathbf{C}^{(1)} & \mathbf{A} + \mathbf{B}^{(1)} \mathbf{K}^{(1)} - \mathbf{L}^{(1)} \mathbf{C}^{(1)}
    \end{bmatrix}}_{\coloneqq\mathbf{\bar A}} \begin{bmatrix}
         \mathbf{x}_{t} \\ \mathbf{z}_{t}^{(1)}
     \end{bmatrix} + \underbrace{\begin{bmatrix}
        \mathbf{I} & 0 \\
        0 & \mathbf{L}^{(1)}
    \end{bmatrix}}_{\coloneqq\mathbf{\bar F}} \begin{bmatrix}
         \mathbf{w}_t \\ \mathbf{v}_t^{(1)}
     \end{bmatrix}.
\end{align*}
The optimal cost of player $1$ is then
\begin{align}
    J^{(1)*} = \mathrm{Tr} \rlpar{\mathbf{\bar P} \mathbf{\bar W}} =  \mathrm{Tr} \rlpar{\mathbf{\bar\Sigma} \mathbf{\bar Q}},
\end{align}
where $$\mathbf{\bar W} =  \mathbf{\bar F} \begin{bmatrix}
        \mathbf{W} & 0 \\
        0 & \mathbf{V}^{(1)}
    \end{bmatrix} \mathbf{\bar F}^\top, \quad\quad \mathbf{\bar Q} =  \begin{bmatrix}
        \mathbf{Q} & 0 \\
        0 & \rlpar{\mathbf{K}^{(1)}}^\top \mathbf{R}^{(1)} \mathbf{K}^{(1)}
    \end{bmatrix},$$ and $\mathbf{\bar P}$, $\mathbf{\bar\Sigma}$ solve the respective Lyapunov equations
\begin{align}
    \mathbf{\bar P} = \mathbf{\bar A}^\top \mathbf{\bar P} \mathbf{\bar A} + \mathbf{\bar Q}, \quad\quad \mathbf{\bar\Sigma} = \mathbf{\bar A} \mathbf{\bar\Sigma} \mathbf{\bar A}^\top + \mathbf{\bar W}.
\end{align}

\subsection{Best Response from Player 2 (evader/maximizer)}\label{sec:best_response_evader}
Having analyzed player $1$'s optimal strategy without player $2$'s control actions in the previous section, we now turn to finding the best response of player $2$ with player $1$'s strategy \eqref{eqn:u_t_1} fixed. The dynamic system now faced by player $2$ incorporates both the state $\mathbf{x}_t$ and player $1$'s state estimate $\mathbf{z}_t^{(1)}$ ($0$-th order belief). We define an augmented dynamic system for player $2$ as
\begin{equation}\label{eqn:xe}
    \mathbf{X}_{t+1}^{(2)} = \mathbf{\bar A}^{(2)} \mathbf{X}_t^{(2)} + \mathbf{\bar B}^{(2)} \mathbf{u}_t^{(2)} + \mathbf{F}^{(2)} \mathbf{\tilde w}_t^{(2)},
\end{equation}
where $\mathbf{X}_t^{(2)} \coloneqq [\mathbf{x}_t^\top \quad (\mathbf{z}_t^{(1)})^\top]^\top \in \bbR^{2n}$, $\mathbf{\tilde  w}_t^{(2)} = [\mathbf{w}_t^\top \quad (\mathbf{v}_t^{(1)})^\top]^\top \in \bbR^{n + p_1}$,
\begin{equation*}
    \mathbf{\bar{A}}^{(2)} = 
    \begin{bmatrix}
        \mathbf{A} & \mathbf{B}^{(1)} \mathbf{K}^{(1)} \\
        \mathbf{L}^{(1)} \mathbf{C}^{(1)} & \mathbf{A} + \mathbf{B}^{(1)} \mathbf{K}^{(1)} - \mathbf{L}^{(1)} \mathbf{C}^{(1)}
    \end{bmatrix}, \quad 
    \mathbf{\bar B}^{(2)} = 
    \begin{bmatrix}
        \mathbf{B}^{(2)} \\ 
        \mathbf{L}^{(1)} \mathbf{D}^{(2)}
    \end{bmatrix}, \quad \mathbf{F}^{(2)} = 
    \begin{bmatrix}
        \mathbf{I} & 0 \\
        0 & \mathbf{L}^{(1)}
    \end{bmatrix}, \quad 
    \mathbf{\tilde W}^{(2)} = 
    \begin{bmatrix}
        \mathbf{W} & 0 \\
        0 & \mathbf{V}^{(1)}
    \end{bmatrix}.
\end{equation*}
Player $2$'s measurement is reformulated as
\begin{equation}\label{eqn:ye}
    \mathbf{y}_t^{(2)} = \mathbf{\bar C}^{(2)} \mathbf{X}_t^{(2)} + \mathbf{v}_t^{(2)},
\end{equation}
where $\mathbf{\bar C}^{(2)} = [\mathbf{C}^{(2)} \quad \mathbf{D}^{(1)} \mathbf{K}^{(1)}]$. Given player $1$'s optimal strategy \eqref{eqn:u_t_1}, the objective function \eqref{eqn:obj_fun} that player $2$ aims to maximize is reformulated as 
\begin{equation*}
    J^{(2)} = \lim_{T \ra \infty} \frac{1}{T}\mathbf{E}_{\mathbf{X}_0^{(2)}, \mathbf{\tilde w}_t^{(2)}, \mathbf{v}_t^{(2)}} \Biggl[\sum_{t=0}^{T-1} \rlpar{\mathbf{X}_t^{(2)}}^\top \mathbf{\bar Q}^{(2)} \mathbf{X}_t^{(2)} + \rlpar{\mathbf{u}_t^{(2)}}^\top \mathbf{R}^{(2)} \mathbf{u}_t^{(2)}\Biggr],
\end{equation*}
where
\begin{equation*}
    \mathbf{\bar Q}^{(2)} = \begin{bmatrix}
        \mathbf{Q} & 0 \\
        0 & \rlpar{\mathbf{K}^{(1)}}^\top \mathbf{R}^{(1)} \mathbf{K}^{(1)}
    \end{bmatrix}.
\end{equation*}

Player $2$'s best response to player $1$'s optimal strategy \eqref{eqn:u_t_1} involves employing a Kalman filter to estimate the \emph{augmented} state $\mathbf{X}_t^{(2)}$, which requires estimating \emph{both} the system state $\mathbf{x}_t$ and the state estimate $\mathbf{z}_t^{(1)}$ of player $1$. Therefore, player $2$ constructs a $1$-st order belief (its belief about player $1$'s belief about the state). Player $2$'s belief state $\mathbf{z}_t^{(2)} \in \bbR^{2n}$ is obtained through a state estimator of the following form
\begin{equation}
    \mathbf{z}_{t+1}^{(2)} = \mathbf{\bar A}^{(2)} \mathbf{z}_t^{(2)} + \mathbf{\bar B}^{(2)} \mathbf{u}_t^{(2)} + \mathbf{L}^{(2)} \rlpar{\mathbf{y}_t^{(2)} - \mathbf{\bar C}^{(2)} \mathbf{z}_t^{(2)}},
\end{equation}
with 
\begin{equation}\label{eqn:L2}
    \mathbf{L}^{(2)} = \mathbf{\bar A}^{(2)} \mathbf{\Sigma}^{(2)} \rlpar{\mathbf{\bar C}^{(2)}}^\top \rlpar{\mathbf{V}^{(2)} + \mathbf{\bar C}^{(2)} \mathbf{\Sigma}^{(2)} \rlpar{\mathbf{\bar C}^{(2)}}^\top}^{-1},
\end{equation}
where the optimal estimator gain $\mathbf{L}^{(2)}$ is constructed from the steady-state estimation error covariance matrix $\mathbf{\Sigma}^{(2)}$ that solves the following algebraic Riccati equation
\begin{align}\label{eqn:Sigma2}
    \mathbf{\Sigma}^{(2)} = \mathbf{F}^{(2)} \mathbf{\tilde W}^{(2)} \rlpar{\mathbf{F}^{(2)}}^\top + \mathbf{\bar A}^{(2)} \mathbf{\Sigma}^{(2)} \rlpar{\mathbf{\bar A}^{(2)}}^\top - \mathbf{\bar A}^{(2)} \mathbf{\Sigma}^{(2)} \rlpar{\mathbf{\bar C}^{(2)}}^\top
     \rlpar{\mathbf{V}^{(2)} + \mathbf{\bar C}^{(2)} \mathbf{\Sigma}^{(2)} \rlpar{\mathbf{\bar C}^{(2)}}^\top}^{-1} \mathbf{\bar C}^{(2)} \mathbf{\Sigma}^{(2)} \rlpar{\mathbf{\bar A}^{(2)}}^\top.
\end{align}
The steady-state estimation error covariance is given by $\mathbf{\Sigma}^{(2)} = \lim_{t \rightarrow \infty} \mathbf{\Sigma}_t^{(2)}$, where $\mathbf{\Sigma}_t^{(2)} = \mathbf{E}[\mathbf{e}_t^{(2)} (\mathbf{e}_t^{(2)})^\top]$ and the estimation error $\mathbf{e}_t^{(2)} \coloneqq \mathbf{X}_t^{(2)} - \mathbf{z}_t^{(2)}$ has dynamics
\begin{align*}
    \mathbf{e}_{t+1}^{(2)} &= \mathbf{X}_{t+1}^{(2)} - \mathbf{z}_{t+1}^{(2)} \\
              &= \rlpar{\mathbf{\bar A}^{(2)} - \mathbf{L}^{(2)} \mathbf{\bar C}^{(2)}} \mathbf{e}_t^{(2)} + \mathbf{F}^{(2)} \mathbf{\tilde w}_t^{(2)} - \mathbf{L}^{(2)} \mathbf{v}_t^{(2)}.
\end{align*}

The best response of player $2$ is the linear (augmented/belief) state feedback
\begin{equation}\label{eqn:u2_t}
    \mathbf{u}_t^{(2)} = \mathbf{K}^{(2)} \mathbf{z}_t^{(2)},
\end{equation}
where the gain matrix $\mathbf{K}^{(2)} \in \bbR^{m_2 \times 2n}$ is given by 
\begin{equation}\label{eqn:K2}
    \mathbf{K}^{(2)} = -\rlpar{\mathbf{R}^{(2)} + \rlpar{\mathbf{\bar B}^{(2)}}^\top \mathbf{P}^{(2)} \mathbf{\bar B}^{(2)}}^{-1} \rlpar{\mathbf{\bar B}^{(2)}}^\top \mathbf{P}^{(2)} \mathbf{\bar A}^{(2)},
\end{equation}
where the cost matrix $\mathbf{P}^{(2)}$ solves the algebraic Riccati equation
\begin{equation}\label{eqn:P2}
    \mathbf{P}^{(2)} = \mathbf{\bar Q}^{(2)} + \rlpar{\mathbf{\bar A}^{(2)}}^\top \mathbf{P}^{(2)} \mathbf{\bar A}^{(2)} - \rlpar{\mathbf{\bar A}^{(2)}}^\top \mathbf{P}^{(2)} \mathbf{\bar B}^{(2)} \rlpar{\mathbf{R}^{(2)} + \rlpar{\mathbf{\bar B}^{(2)}}^\top \mathbf{P}^{(2)} \mathbf{\bar B}^{(2)}}^{-1} \rlpar{\mathbf{\bar B}^{(2)}}^\top \mathbf{P}^{(2)} \mathbf{\bar A}^{(2)}.
\end{equation}
There exist unique solutions $\mathbf{\Sigma}^{(2)}$ and $\mathbf{P}^{(2)}$ of the Riccati equations if the \emph{augmented} system represented by \eqref{eqn:xe} and \eqref{eqn:ye} is stabilizable and detectable. Moreover, we emphasize that player $2$'s optimal cost is bounded only when a solution exists where $\mathbf{R}^{(2)} + (\mathbf{\bar B}^{(2)})^\top \mathbf{P}^{(2)} \mathbf{\bar B}^{(2)} \prec 0$, which requires that the penalty on player $2$'s control actions is sufficiently negative definite.

The closed-loop augmented state and belief state dynamics of player $2$ are formulated as
\begin{align*} \nonumber
     \begin{bmatrix}
         \mathbf{X}^{(2)}_{t+1} \\ \mathbf{z}_{t+1}^{(2)}
     \end{bmatrix} = \underbrace{\begin{bmatrix}
        \mathbf{\bar A}^{(2)} & \mathbf{\bar B}^{(2)} \mathbf{K}^{(2)} \\
        \mathbf{L}^{(2)} \mathbf{\bar C}^{(2)} & \mathbf{\bar A}^{(2)} + \mathbf{\bar B}^{(2)} \mathbf{K}^{(2)} - \mathbf{L}^{(2)} \mathbf{\bar C}^{(2)}
    \end{bmatrix}}_{\coloneqq\mathbf{\tilde A}} \begin{bmatrix}
         \mathbf{X}^{(2)}_{t} \\ \mathbf{z}_{t}^{(2)}
     \end{bmatrix} + \underbrace{\begin{bmatrix}
        \mathbf{F}^{(2)} & 0 \\
        0 & \mathbf{L}^{(2)}
    \end{bmatrix}}_{\coloneqq\mathbf{\tilde F}} \begin{bmatrix}
         \mathbf{\tilde w}_t^{(2)} \\ \mathbf{v}_t^{(2)}
     \end{bmatrix}.
\end{align*}
The optimal cost of player $2$ is then
\begin{align}
    J^{(2)*} = \mathrm{Tr} \rlpar{\mathbf{\tilde P} \mathbf{\tilde W}} = \mathrm{Tr} \rlpar{\mathbf{\tilde\Sigma} \mathbf{\tilde Q}},
\end{align}
where $$\mathbf{\tilde W} = \mathbf{\tilde F} \begin{bmatrix}
        \mathbf{\bar W}^{(2)} & 0 \\
        0 & \mathbf{V}^{(2)}
    \end{bmatrix} \mathbf{\tilde F}^\top, \quad\quad \mathbf{\tilde Q} =  \begin{bmatrix}
        \mathbf{\bar Q}^{(2)} & 0 \\
        0 & \rlpar{\mathbf{K}^{(2)}}^\top \mathbf{R}^{(2)} \mathbf{K}^{(2)}
    \end{bmatrix},$$ and $\mathbf{\tilde P}$, $\mathbf{\tilde\Sigma}$ solve the respective Lyapunov equations
\begin{align}
    \mathbf{\tilde P} = \mathbf{\tilde A}^\top \mathbf{\tilde P} \mathbf{\tilde A} + \mathbf{\tilde Q}, \quad\quad \mathbf{\tilde\Sigma} = \mathbf{\tilde A} \mathbf{\tilde\Sigma} \mathbf{\tilde A}^\top + \mathbf{\tilde W}. 
\end{align}

\subsection{General Best Response Dynamics}
Based on the analysis in Section \ref{sec:lqg_pursuer} and \ref{sec:best_response_evader}, we now derive explicit expressions for the general best response dynamics of each player that, upon convergence, yield a Nash equilibrium within our pure linear dynamic output feedback strategy class. Below we summarize the best response dynamics for each player for our stochastic dynamic game with partial and asymmetric information. In each iteration, the opposing player's strategy is fixed, and each player faces an LQG control problem, but for an augmented system that includes both the system state and the opposing player's increasing higher-order belief states. Each iteration increases the dimension of the state estimate for the augmented belief state system to an incremental integer multiple of the original state dimension. This reformulation captures recursive definitions for the augmented belief state dynamics, objective function, and state estimator. This process iterates until both players' best responses converge. If the iterations converge, we obtain a Nash equilibrium within our strategy class \cite{bacsar1998dynamic}, that is, no player can unilaterally improve their cost by deviating from their equilibrium policy.

The expressions for the best response dynamics below follow the recursive solving of an LQG problem for each player. By construction, augmented belief dynamics considered by each player increase in dimension at each best response iteration, thus increasing the internal dimensions of the feedback strategies employed by each player.

\subsubsection{Player 1's Best Response Characteristics}\label{sec:general_best_response_pursuer}
Player $1$'s best response at iteration $k = 2,3,\ldots$ considers the augmented belief dynamics and measurement 
\begin{equation}\label{eqn:pursuer_sys_k}
    \begin{split}
        \mathbf{X}_{t+1,k}^{(1)} &= \mathbf{\bar A}_k^{(1)} \mathbf{X}_{t,k}^{(1)} + \mathbf{\bar B}_k^{(1)} \mathbf{u}_{t,k}^{(1)} + \mathbf{F}_k^{(1)} \mathbf{\tilde{w}}_t^{(1)}, \\
        \mathbf{y}_{t,k}^{(1)} &= \mathbf{\bar C}_k^{(1)} \mathbf{X}_{t,k}^{(1)} + \mathbf{v}_t^{(1)},
    \end{split}
\end{equation}
where $\mathbf{X}_{t,k}^{(1)} \coloneqq [\mathbf{x}_{t}^\top \quad (\mathbf{z}_{t,k-1}^{(2)})^\top]^\top \in \mathbb{R}^{(2k-1)n}$, $\mathbf{\tilde{w}}_{t}^{(1)} = [\mathbf{w}_t^\top \quad (\mathbf{v}_t^{(2)})^\top]^\top \in \bbR^{n + p_2}$,
\begin{align*}
    \mathbf{\bar A}_{k}^{(1)} =
    \begin{bmatrix}
        \mathbf{A} & \mathbf{B}^{(2)} \mathbf{K}_{k-1}^{(2)} \\
        \mathbf{L}_{k-1}^{(2)} \mathbf{C}^{(2)} & \mathbf{\bar A}_{k-1}^{(2)} + \mathbf{\bar B}_{k-1}^{(2)} \mathbf{K}_{k-1}^{(2)} - \mathbf{L}_{k-1}^{(2)} \mathbf{\bar C}^{(2)}_{k-1}
    \end{bmatrix}, \quad \mathbf{\bar B}_k^{(1)} =
    \begin{bmatrix}
        \mathbf{B}^{(1)} \\
        \mathbf{L}_{k-1}^{(2)} \mathbf{D}^{(1)}
    \end{bmatrix}, \quad
    \mathbf{\bar C}_k^{(1)} = 
    \begin{bmatrix}
        \mathbf{C}^{(1)} & \mathbf{D}^{(2)} \mathbf{K}_{k-1}^{(2)}
    \end{bmatrix}, 
\end{align*}
\begin{align*}
    \mathbf{F}_{k}^{(1)} = 
    \begin{bmatrix}
        \mathbf{I} & 0 \\
        0 & \mathbf{L}_{k-1}^{(2)}
    \end{bmatrix}, \quad
    \mathbf{\tilde W}^{(1)} =
    \begin{bmatrix}
        \mathbf{W} & 0 \\
        0 & \mathbf{V}^{(2)} 
    \end{bmatrix}.
\end{align*}
These matrices' dimensions at iteration $k$ are $\mathbf{\bar A}_k^{(1)} \in \bbR^{(2k-1)n \times (2k-1)n}$, $\mathbf{\bar B}_k^{(1)} \in \bbR^{(2k-1)n \times m_1}$, $\mathbf{\bar C}_k^{(1)} \in \bbR^{p_1 \times (2k-1)n}$, $\mathbf{F}_k^{(1)} \in \bbR^{(2k-1)n \times n+p_1}$, and $\mathbf{\tilde W}^{(1)} \in \bbR^{n+p_2 \times n+p_2}$. The parameters $\mathbf{K}_{k-1}^{(2)}$, $\mathbf{L}_{k-1}^{(2)}$, $\mathbf{\bar A}_{k-1}^{(2)}$, $\mathbf{\bar B}_{k-1}^{(2)}$, and $\mathbf{\bar C}_{k-1}^{(2)}$ are described in Section \ref{sec:general_best_response_evader} with $\mathbf{K}_1^{(2)} = \mathbf{K}^{(2)}$, $\mathbf{L}_1^{(2)} = \mathbf{L}^{(2)}$, $\mathbf{\bar A}_1^{(2)} = \mathbf{\bar A}^{(2)}$, $\mathbf{\bar B}_1^{(2)} = \mathbf{\bar B}^{(2)}$, and $\mathbf{\bar C}_1^{(2)} = \mathbf{\bar C}^{(2)}$ from player $2$'s first best response in Section \ref{sec:best_response_evader}.

The objective function for player $1$ to minimize at iteration $k$ is
\begin{equation} \nonumber
    J_k^{(1)} = \lim_{T \ra \infty} \frac{1}{T} \mathbf{E}_{\mathbf{X}_0^{(1)}, \mathbf{\tilde w}_t^{(1)}, \mathbf{v}_t^{(1)}} \left[\sum_{t=0}^{T-1} \rlpar{\mathbf{X}_{t,k}^{(1)}}^\top \mathbf{\bar Q}_k^{(1)} \mathbf{X}_{t,k}^{(1)}  + \rlpar{\mathbf{u}_{t,k}^{(1)}}^\top \mathbf{R}^{(1)} \mathbf{u}_{t,k}^{(1)} \right],
\end{equation}
where
\begin{align*}
    \mathbf{\bar Q}_k^{(1)} = 
    \begin{bmatrix}
        \mathbf{Q} & 0 \\
        0 & \rlpar{\mathbf{K}_{k-1}^{(2)}}^\top \mathbf{R}^{(2)} \mathbf{K}_{k-1}^{(2)}
    \end{bmatrix} \in \bbR^{(2k-1)n \times (2k-1)n}.
\end{align*}
Player $1$'s belief state $\mathbf{z}_{t,k}^{(1)} \in \mathbb{R}^{(2k-1)n}$ at iteration $k$ is obtained through a state estimator
\begin{align}
    \mathbf{z}_{t+1,k}^{(1)} = \mathbf{\bar A}_k^{(1)} \mathbf{z}_{t,k}^{(1)} + \mathbf{\bar B}_k^{(1)} \mathbf{u}_{t,k}^{(1)} + \mathbf{L}_{k}^{(1)} \rlpar{\mathbf{y}_{t,k}^{(1)} - \mathbf{\bar C}_k^{(1)} \mathbf{z}_{t,k}^{(1)}},
\end{align}
and the best response of player $1$ at iteration $k$ is the linear (augmented/belief) state feedback
\begin{equation}\label{eqn:u_tk^p}
    \mathbf{u}_{t,k}^{(1)} = \mathbf{K}_k^{(1)} \mathbf{z}_{t,k}^{(1)},
\end{equation}
where the optimal feedback gain $\mathbf{K}_k^{(1)} \in \bbR^{m_1 \times (2k-1)n}$, the state estimator gain $\mathbf{L}_k^{(1)} \in  \bbR^{(2k-1)n \times p_1}$, the associated cost matrix $\mathbf{P}_k^{(1)} \in \bbR^{(2k-1)n \times (2k-1)n}$, and the estimation error covariance matrix $\mathbf{\Sigma}_k^{(1)} \in \bbR^{(2k-1)n \times (2k-1)n}$ are given by
\begin{align}
    \mathbf{K}_k^{(1)} &= -\rlpar{\mathbf{R}^{(1)} + \rlpar{\mathbf{\bar B}_k^{(1)}}^\top \mathbf{P}_k^{(1)} \mathbf{\bar B}_k^{(1)}}^{-1} \rlpar{\mathbf{\bar B}_k^{(1)}}^\top \mathbf{P}_k^{(1)} \mathbf{\bar A}_k^{(1)}, \label{eqn:K1_k} \\
    \mathbf{P}_k^{(1)} &= \mathbf{\bar Q}_k^{(1)} + \rlpar{\mathbf{\bar A}_k^{(1)}}^\top \mathbf{P}_k^{(1)} \mathbf{\bar A}_k^{(1)} - \rlpar{\mathbf{\bar A}_k^{(1)}}^\top \mathbf{P}_k^{(1)} \mathbf{\bar B}_k^{(1)} \rlpar{\mathbf{R}^{(1)} + \rlpar{\mathbf{\bar B}_k^{(1)}}^\top \mathbf{P}_k^{(1)} \mathbf{\bar B}_k^{(1)}}^{-1} \rlpar{\mathbf{\bar B}_k^{(1)}}^\top \mathbf{P}_k^{(1)} \mathbf{\bar A}_k^{(1)}, \label{eqn:P1_k} \\
    \mathbf{L}_k^{(1)} &= \mathbf{\bar A}_k^{(1)} \mathbf{\Sigma}_k^{(1)} \rlpar{\mathbf{\bar C}_k^{(1)}}^\top \rlpar{\mathbf{V}^{(1)} + \mathbf{\bar C}_k^{(1)} \mathbf{\Sigma}_k^{(1)} \rlpar{\mathbf{\bar C}_k^{(1)}}^\top}^{-1}, \label{eqn:L1_k}\\
    \mathbf{\Sigma}_k^{(1)} &= \mathbf{F}_k^{(1)} \mathbf{\tilde W}^{(1)} \rlpar{\mathbf{F}_k^{(1)}}^\top + \mathbf{\bar A}_k^{(1)} \mathbf{\Sigma}_k^{(1)} \rlpar{\mathbf{\bar A}_k^{(1)}}^\top - \mathbf{\bar A}_k^{(1)} \mathbf{\Sigma}_k^{(1)} \rlpar{\mathbf{\bar C}_k^{(1)}}^\top \rlpar{\mathbf{V}^{(1)} + \mathbf{\bar C}_k^{(1)} \mathbf{\Sigma}_k^{(1)} \rlpar{\mathbf{\bar C}_k^{(1)}}^\top}^{-1} \mathbf{\bar C}_k^{(1)} \mathbf{\Sigma}_k^{(1)} \rlpar{\mathbf{\bar A}_k^{(1)}}^\top. \label{eqn:Sigma1_k}
\end{align}

Although the penalty matrix $\mathbf{\bar Q}_k^{(1)}$ is indefinite due to $\mathbf{Q} \succeq 0$ and $(\mathbf{K}_{k-1}^{(2)})^\top \mathbf{R}^{(2)} \mathbf{K}_{k-1}^{(2)} \prec 0$, unique solutions $\mathbf{P}_k^{(1)}$ and $\mathbf{\Sigma}_k^{(1)}$ for the Riccati equations exist under certain conditions related to stabilizability of $(\bar{\mathbf{A}}_k^{(1)},\ \bar{\mathbf{B}}_k^{(1)})$, $\mathbf{R}^{(1)} \succ 0$ and detectability of $(\bar{\mathbf{A}}_k^{(1)},\ \bar{\mathbf{C}}_k^{(1)})$, $\mathbf{V}^{(1)} \succ 0$, and how they relate to the indefiniteness of $\mathbf{\bar Q}_k^{(1)}$ \cite{willems2003least,molinari1973stable}. Existence of solutions is verified numerically. The closed-loop augmented state and belief state dynamics of player $1$ at iteration $k$ are
\begin{align*} \nonumber
     \begin{bmatrix}
         \mathbf{X}^{(1)}_{t+1,k} \\ \mathbf{z}_{t+1,k}^{(1)}
     \end{bmatrix} = \underbrace{\begin{bmatrix}
        \mathbf{\bar A}_k^{(1)} & \mathbf{\bar B}_k^{(1)} \mathbf{K}_k^{(1)} \\
        \mathbf{L}_k^{(1)} \mathbf{\bar C}_k^{(1)} & \mathbf{\bar A}_k^{(1)} + \mathbf{\bar B}_k^{(1)} \mathbf{K}_k^{(1)} - \mathbf{L}_k^{(1)} \mathbf{\bar C}_k^{(1)}
    \end{bmatrix}}_{\coloneqq\mathbf{\bar A}_k} \begin{bmatrix}
         \mathbf{X}^{(1)}_{t,k} \\ \mathbf{z}_{t,k}^{(1)}
     \end{bmatrix} + \underbrace{\begin{bmatrix}
        \mathbf{F}_k^{(1)} & 0 \\
        0 & \mathbf{L}_k^{(1)}
    \end{bmatrix}}_{\coloneqq\mathbf{\bar F}_k} \begin{bmatrix}
         \mathbf{\tilde w}_t^{(1)} \\ \mathbf{v}_t^{(1)}
     \end{bmatrix}.
\end{align*}
The optimal cost of player $1$ at iteration $k$ is then
\begin{align}
    J_k^{(1)*} = \mathrm{Tr} \rlpar{\mathbf{\bar P}_k \mathbf{\bar W}_k} =  \mathrm{Tr} \rlpar{\mathbf{\bar\Sigma}_k \mathbf{\bar Q}_k},
\end{align}
where $$\mathbf{\bar W}_k = \mathbf{\bar F}_k \begin{bmatrix}
    \mathbf{\tilde W}^{(1)} & 0 \\
    0 & \mathbf{V}^{(1)}
\end{bmatrix} \mathbf{\bar F}_k^\top, \quad\quad \mathbf{\bar Q}_k =  \begin{bmatrix}
    \mathbf{\bar Q}_k^{(1)} & 0 \\
    0 & \rlpar{\mathbf{K}_k^{(1)}}^\top \mathbf{R}^{(1)} \mathbf{K}_k^{(1)}
\end{bmatrix},$$ and $\mathbf{\bar P}_k$, $\mathbf{\bar\Sigma}_k$ solve the respective Lyapunov equations
\begin{align}
    \mathbf{\bar P}_k = \mathbf{\bar A}_k^\top \mathbf{\bar P}_k \mathbf{\bar A}_k + \mathbf{\bar Q}_k, \quad\quad \mathbf{\bar\Sigma}_k = \mathbf{\bar A}_k \mathbf{\bar\Sigma}_k \mathbf{\bar A}_k^\top + \mathbf{\bar W}_k. 
\end{align}

\subsubsection{Player 2's Best Response Characteristics}\label{sec:general_best_response_evader}
Player $2$'s best response at iteration $k=1,2,\ldots$ considers the augmented belief dynamics and measurement
\begin{align}\label{eqn:max_aug_sys}
    \begin{split}
        \mathbf{X}_{t+1,k}^{(2)} &= \mathbf{\bar A}_k^{(2)} \mathbf{X}_{t,k}^{(2)} + \mathbf{\bar B}_k^{(2)} \mathbf{u}_{t,k}^{(2)} + \mathbf{F}_k^{(2)} \mathbf{\tilde w}_t^{(2)}, \\
        \mathbf{y}_{t,k}^{(2)} &= \mathbf{\bar C}_k^{(2)} \mathbf{X}_{t,k}^{(2)} + \mathbf{v}_t^{(2)},
    \end{split}
\end{align}
where $\mathbf{X}_{t,k}^{(2)} \coloneqq [\mathbf{x}_{t}^\top \quad (\mathbf{z}_{t,k}^{(1)})^\top]^\top \in \mathbb{R}^{2kn}$, $\mathbf{\tilde{w}}_{t}^{(2)} = [\mathbf{w}_t^\top \quad (\mathbf{v}_t^{(1)})^\top]^\top \in \bbR^{n+p_1}$,
\begin{align*}
    \mathbf{\bar A}_{k}^{(2)} =
    \begin{bmatrix}
        \mathbf{A} & \mathbf{B}^{(1)} \mathbf{K}_k^{(1)} \\
        \mathbf{L}_{k}^{(1)} \mathbf{C}^{(1)} & \mathbf{\bar A}_{k}^{(1)} + \mathbf{\bar B}_{k}^{(1)} \mathbf{K}_{k}^{(1)} - \mathbf{L}_{k}^{(1)} \mathbf{\bar C}_{k}^{(1)}
    \end{bmatrix}, \quad
    \mathbf{\bar B}_k^{(2)} =
    \begin{bmatrix}
        \mathbf{B}^{(2)} \\
        \mathbf{L}_k^{(1)} \mathbf{D}^{(2)}
    \end{bmatrix}, \quad \mathbf{\bar C}_k^{(2)} = 
    \begin{bmatrix}
        \mathbf{C}^{(2)} & \mathbf{D}^{(1)} \mathbf{K}_k^{(1)}
    \end{bmatrix},
\end{align*}
\begin{align*}
    \mathbf{F}_{k}^{(2)} = 
    \begin{bmatrix}
        \mathbf{I} & 0 \\
        0 & \mathbf{L}_k^{(1)}
    \end{bmatrix}, \quad
    \mathbf{\tilde W}^{(2)} =
    \begin{bmatrix}
        \mathbf{W} & 0 \\
        0 & \mathbf{V}^{(1)} 
    \end{bmatrix}.
\end{align*}
These matrices' dimensions at iteration $k$ are $\mathbf{\bar A}_k^{(2)} \in \bbR^{2kn \times 2kn}$, $\mathbf{\bar B}_k^{(2)} \in \bbR^{2kn \times m_2}$, $\mathbf{\bar C}_k^{(2)} \in \bbR^{p_2 \times 2kn}$, $\mathbf{F}_k^{(2)} \in \bbR^{2kn \times n+p_1}$, and $\mathbf{\tilde W}^{(2)} \in \bbR^{n+p_1 \times n+p_1}$. The parameters $\mathbf{K}_k^{(1)}$, $\mathbf{L}_k^{(1)}$, $\mathbf{\bar A}_k^{(1)}$, $\mathbf{\bar B}_k^{(1)}$, and $\mathbf{\bar C}_k^{(1)}$ are presented in Section \ref{sec:general_best_response_pursuer} with $\mathbf{K}_1^{(1)} = \mathbf{K}^{(1)}$, $\mathbf{L}_1^{(1)} = \mathbf{L}^{(1)}$, $\mathbf{\bar A}_1^{(1)} = \mathbf{A}$, $\mathbf{\bar B}_1^{(1)} = \mathbf{B}^{(1)}$, and $\mathbf{C}_1^{(1)} = \mathbf{C}^{(1)}$ from player $1$'s optimal strategy in Section \ref{sec:lqg_pursuer}.

The objective function for player $2$ to maximize at iteration $k$ is 
\begin{equation*}
    J_k^{(2)} = \lim_{T \ra \infty} \frac{1}{T} \mathbf{E}_{\mathbf{X}_0^{(2)}, \mathbf{\tilde w}_t^{(2)}, \mathbf{v}_t^{(2)}} \left[\sum_{t=0}^{T-1} \rlpar{\mathbf{X}_{t,k}^{(2)}}^\top \mathbf{\bar Q}_k^{(2)} \mathbf{X}_{t,k}^{(2)} + \rlpar{\mathbf{u}_{t,k}^{(2)}}^\top \mathbf{R}^{(2)} \mathbf{u}_{t,k}^{(2)} \right],
\end{equation*}
where
\begin{align*}
    \mathbf{\bar Q}_k^{(2)} = 
    \begin{bmatrix}
        \mathbf{Q} & 0 \\
        0 & \rlpar{\mathbf{K}_k^{(1)}}^\top \mathbf{R}^{(1)} \mathbf{K}_k^{(1)}
    \end{bmatrix} \in \bbR^{2kn \times 2kn}.
\end{align*}
Player $2$'s belief state $\mathbf{z}_{t,k}^{(2)} \in \mathbb{R}^{2kn}$ at iteration $k$ is obtained through a state estimator
\begin{align}
    \mathbf{z}_{t+1,k}^{(2)} = \mathbf{\bar A}_k^{(2)} \mathbf{z}_{t,k}^{(2)} + \mathbf{\bar B}_k^{(2)} \mathbf{u}_{t,k}^{(2)} + \mathbf{L}_{k}^{(2)} \rlpar{\mathbf{y}_{t,k}^{(2)} - \mathbf{\bar C}_k^{(2)} \mathbf{z}_{t,k}^{(2)}},
\end{align}
and the best response of player $2$ at iteration $k$ is the linear (augmented/belief) state feedback
\begin{equation}
    \mathbf{u}_{t,k}^{(2)} = \mathbf{K}_k^{(2)} \mathbf{z}_{t,k}^{(2)},
\end{equation}
 where the optimal feedback gain $\mathbf{K}_k^{(2)} \in \bbR^{m_2 \times 2kn}$, the state estimator gain $\mathbf{L}_k^{(2)} \in \bbR^{2kn \times p_2}$, the associated cost matrix $\mathbf{P}_k^{(2)} \in \bbR^{2kn \times 2kn}$, and the estimation error covariance matrix $\mathbf{\Sigma}_k^{(2)} \in \bbR^{2kn \times 2kn}$ are given by
 \begin{align}
    \mathbf{K}_k^{(2)} &= -\rlpar{\mathbf{R}^{(2)} + \rlpar{\mathbf{\bar B}_k^{(2)}}^\top \mathbf{P}_k^{(2)} \mathbf{\bar B}_k^{(2)}}^{-1} \rlpar{\mathbf{\bar B}_k^{(2)}}^\top \mathbf{P}_k^{(2)} \mathbf{\bar A}_k^{(2)}, \label{eqn:K2_k} \\
    \mathbf{P}_k^{(2)} &= \mathbf{\bar Q}_k^{(2)} + \rlpar{\mathbf{\bar A}_k^{(2)}}^\top \mathbf{P}_k^{(2)} \mathbf{\bar A}_k^{(2)} - \rlpar{\mathbf{\bar A}_k^{(2)}}^\top \mathbf{P}_k^{(2)} \mathbf{\bar B}_k^{(2)} \rlpar{\mathbf{R}^{(2)} + \rlpar{\mathbf{\bar B}_k^{(2)}}^\top \mathbf{P}_k^{(2)} \mathbf{\bar B}_k^{(2)}}^{-1} \rlpar{\mathbf{\bar B}_k^{(2)}}^\top \mathbf{P}_k^{(2)} \mathbf{\bar A}_k^{(2)}, \label{eqn:P2_k} \\
    \mathbf{L}_k^{(2)} &= \mathbf{\bar A}_k^{(2)} \mathbf{\Sigma}_k^{(2)} \rlpar{\mathbf{\bar C}_k^{(2)}}^\top \rlpar{\mathbf{V}^{(2)} + \mathbf{\bar C}_k^{(2)} \mathbf{\Sigma}_k^{(2)} \rlpar{\mathbf{\bar C}_k^{(2)}}^\top}^{-1}, \label{eqn:L2_k}\\
    \mathbf{\Sigma}_k^{(2)} &= \mathbf{F}_k^{(2)} \mathbf{\tilde W}^{(2)} \rlpar{\mathbf{F}_k^{(2)}}^\top + \mathbf{\bar A}_k^{(2)} \mathbf{\Sigma}_k^{(2)} \rlpar{\mathbf{\bar A}_k^{(2)}}^\top - \mathbf{\bar A}_k^{(2)} \mathbf{\Sigma}_k^{(2)} \rlpar{\mathbf{\bar C}_k^{(2)}}^\top \rlpar{\mathbf{V}^{(2)} + \mathbf{\bar C}_k^{(2)} \mathbf{\Sigma}_k^{(2)} \rlpar{\mathbf{\bar C}_k^{(2)}}^\top}^{-1} \mathbf{\bar C}_k^{(2)} \mathbf{\Sigma}_k^{(2)} \rlpar{\mathbf{\bar A}_k^{(2)}}^\top. \label{eqn:Sigma2_k}
\end{align}
For player $2$'s augmented LQG problem, the solvability of the Riccati equation corresponds to $\mathbf{P}_k^{(2)}$ requires $\mathbf{R}^{(2)} + (\mathbf{\bar B}_k^{(2)})^\top \mathbf{P}_k^{(2)} \mathbf{\bar B}_k^{(2)} \prec 0$, ensuring the game's upper value remains bounded. This condition, which requires $\mathbf{R}^{(2)}$ to be sufficiently negative definite is established in \cite{bacsar1998dynamic} and is verified numerically at each best response iteration.

The closed-loop augmented state and belief state dynamics of player $2$ at iteration $k$ are
\begin{align*} \nonumber
     \begin{bmatrix}
         \mathbf{X}^{(2)}_{t+1,k} \\ \mathbf{z}_{t+1,k}^{(2)}
     \end{bmatrix} = \underbrace{\begin{bmatrix}
        \mathbf{\bar A}_k^{(2)} & \mathbf{\bar B}_k^{(2)} \mathbf{K}_k^{(2)} \\
        \mathbf{L}_k^{(2)} \mathbf{\bar C}_k^{(2)} & \mathbf{\bar A}_k^{(2)} + \mathbf{\bar B}_k^{(2)} \mathbf{K}_k^{(2)} - \mathbf{L}_k^{(2)} \mathbf{\bar C}_k^{(2)}
    \end{bmatrix}}_{\coloneqq\mathbf{\tilde A}_k} \begin{bmatrix}
         \mathbf{X}^{(2)}_{t,k} \\ \mathbf{z}_{t,k}^{(2)}
     \end{bmatrix} + \underbrace{\begin{bmatrix}
        \mathbf{F}_k^{(2)} & 0 \\
        0 & \mathbf{L}_k^{(2)}
    \end{bmatrix}}_{\coloneqq\mathbf{\tilde F}_k} \begin{bmatrix}
          \mathbf{\tilde w}_t^{(2)} \\ \mathbf{v}_t^{(2)}
     \end{bmatrix}.
\end{align*}
The optimal cost of player $2$ at iteration $k$ is then
\begin{align}
    J_k^{(2)*} = \mathrm{Tr} \rlpar{\mathbf{\tilde P}_k \mathbf{\tilde W}_k} =  \mathrm{Tr} \rlpar{\mathbf{\tilde\Sigma}_k \mathbf{\tilde Q}_k},
\end{align}
where $$\mathbf{\tilde W}_k = \mathbf{\tilde F}_k \begin{bmatrix}
    \mathbf{\tilde W}^{(2)} & 0 \\
    0 & \mathbf{V}^{(2)}
\end{bmatrix} \mathbf{\tilde F}_k^\top, \quad\quad \mathbf{\tilde Q}_k =  \begin{bmatrix}
    \mathbf{\bar Q}_k^{(2)} & 0 \\
    0 & \rlpar{\mathbf{K}_k^{(2)}}^\top \mathbf{R}^{(2)} \mathbf{K}_k^{(2)}
\end{bmatrix},$$ and $\mathbf{\tilde P}_k$, $\mathbf{\tilde\Sigma}_k$ solve the respective Lyapunov equations
\begin{align}
    \mathbf{\tilde P}_k = \mathbf{\tilde A}_k^\top \mathbf{\tilde P}_k \mathbf{\tilde A}_k + \mathbf{\tilde Q}_k, \quad\quad \mathbf{\tilde\Sigma}_k = \mathbf{\tilde A}_k \mathbf{\tilde\Sigma}_k \mathbf{\tilde A}_k^\top + \mathbf{\tilde W}_k. 
\end{align}

In theory, the internal state dimensions of both players' feedback strategies grow towards infinity as the best response dynamics evolve. However, our numerical results in the next section show that the game's value converges rapidly within a few iterations. To help explain this convergence, our previous analysis \cite{guan2025best} shows that the controllability and observability Gramian eigenvalues and Hankel singular values of both players' higher-order belief dynamics decay rapidly. Theorem~1 in \cite{guan2025best}, summarized in the Appendix, bounds these decay rates using Cholesky estimates and shows that higher-order belief dynamics can be approximated by low-order ones with bounded error. As a result, equilibrium strategies based on low-order belief states give a good approximation of the true Nash equilibrium strategies. This result was established and numerically validated for open-loop stable systems without feed-through.

\begin{remark}
    Theorem~1 is not limited to open-loop stable systems without feed-through. Theorem~1's conditions, namely stability, diagonalizability, controllability, and observability of the augmented belief dynamics matrix at each iteration $k$, are stated on the iteration-$k$ matrices $(\mathbf{\bar{A}}_k^{(i)}, \mathbf{\bar{B}}_k^{(i)}, \mathbf{\bar{C}}_k^{(i)})$ themselves and do not directly depend on the feed-through matrices $\mathbf{D}^{(i)}$ or on the open-loop stability of $A$. Theorem~1 therefore applies whether or not $\mathbf{D}^{(i)}=0$, as long as these conditions hold for the resulting matrices. Regarding open-loop stability, the only iteration where this becomes a binding constraint is $k=1$, since player~1's augmented belief matrix at this iteration reduces to $A$ itself, prior to any feedback being incorporated. If $A$ is open-loop unstable, Theorem~1 does not apply at iteration $k=1$, because the Gramian equations in Theorem~1 require the underlying system to be Schur stable \cite{antoulas2005approximation}. At later iterations, however, the augmented belief dynamics matrix incorporates the feedback and estimator gains from prior best responses and generally differs from $A$, so Theorem~1 applies at any iteration $k\geq2$ for which its stated conditions hold, regardless of whether $A$ itself was open-loop stable.  
\end{remark}

\begin{remark}
    A general convergence proof for the best response iteration remains open. The main difficulty is that the space of output-feedback strategies under asymmetric information has no known finite characterization, so standard fixed-point arguments do not directly apply. Two analytical facts still constrain the iteration's behavior. First, each player's own subproblem cost is monotone at every iteration (non-decreasing for the pursuer, non-increasing for the evader), since each best response is an exact LQG optimum given the opponent's fixed strategy. This follows from construction, not simulation. Second, Theorem 1 (in Appendix) gives an explicit bound showing that higher-order belief dynamics can be approximated by low-order ones with bounded error, which explains analytically why additional belief order yields diminishing returns. These results support, but do not prove, the convergence observed in our experiments. Establishing sufficient conditions on the system and penalty matrices under which convergence is guaranteed remains a challenging open problem for future work.
\end{remark}

\section{Numerical Experiments}\label{sec:experimental_results}
In this section, we analyze how various factors affect the average infinite horizon cost in a planar pursuit-evasion game under partial and asymmetric information settings, using dynamical and sensor models of Section \ref{sec:two_player_pe_game}. Specifically, we examine the effects of players' asymmetric beliefs, belief orders, relative controllability and observability, and ECM attack on the game's value. We begin by presenting the convergence results of both players' best response strategies in Section \ref{sec:optimal_cost}. Next, we investigate how asymmetric beliefs affect the game cost and the players' mean trajectories in Section \ref{sec:perf_vs_partial}. We then explore the impact of belief orders and relative controllability and observability on both the game cost and the players' trajectories in Sections \ref{sec:belief_order_impact} and \ref{sec:pursuit_evasion_maneuver_observe}. Finally, we assess how the evader’s ECM attack on the pursuer’s sensor affects the game’s cost in Section \ref{sec:pursuit_evasion_game_attack}. Sections IV.A-D set $\mathbf{D}^{(p)} = \mathbf{{D}^{(e)}} = 0$, corresponding to no sensor interference; Section IV.E introduces nonzero feed-through matrices to study the strategic impact of ECM attacks.

We consider the in-plane pursuit-evasion game described in Section \ref{sec:two_player_pe_game}, restricting attention to the $x$-$y$ plane of the LVLH frame. Setting $n=4$, $\Delta t = 0.1 \mathrm{s}$, and $w\approx1.13\times10^{-3} /\mathrm{s}$, the in-plane state $\mathbf{x}_t = [x_t,\ y_t,\ \dot{x}_t,\ \dot{y}_t]^\top\in\mathbb{R}^4$ and the system matrices reduce to
\begin{equation*}
    \mathbf{A} \approx \begin{bmatrix}
        1 & 0 & \Delta t & 0 \\
        0 & 1 & 0 & \Delta t \\
        0 & 0 & 1 & 0 \\
        0 & 0 & 0 & 1
    \end{bmatrix}, \quad  
    \mathbf{B}^{(p)} = \mathbf{G}^{(p)} \approx \begin{bmatrix}
        \frac{1}{2}\Delta t^2 & 0 \\
        0 & \frac{1}{2}\Delta t^2  \\
        \Delta t & 0 \\
        0 & \Delta t
    \end{bmatrix}, \quad \mathbf{B}^{(e)} = -\mathbf{B}^{(p)}, \quad 
    \mathbf{G}^{(e)} = -\mathbf{G}^{(p)}, 
\end{equation*}
\begin{equation*}
    \mathbf{C}^{(p)} = \mathbf{C}^{(e)} = \begin{bmatrix}
        1 & 0 & 0 & 0 \\
        0 & 1 & 0 & 0
    \end{bmatrix}, \quad
    \mathbf{D}^{(p)} = \mathbf{D}^{(e)} = \mathbf{0}_{2},
\end{equation*}
and we define the noise covariances, penalty matrices, and initial conditions as
\begin{equation*}
    \mathbf{W}_d^{(p)} = \mathbf{W}_d^{(e)} = 10^2\mathbf{I}_2, \quad \mathbf{V}^{(p)} = \mathbf{V}^{(e)} = 10^{2} \mathbf{I}_2, \quad \mathbf{Q} = 0.1 \mathbf{I}_4, \quad \mathbf{R}_u^{(p)} = \mathbf{I}_2, \quad \mathbf{R}_u^{(e)} = -10 \mathbf{I}_2,
\end{equation*}
\begin{equation*}
    \mathbf{R}_s^{(p)} = \mathbf{0}_2, \quad \mathbf{R}_s^{(e)} = \mathbf{0}_2, \quad \mathbf{x}_0 = [-100,\ 0,\ 0,\ 0]^\top\mathrm{m}, \quad \mathbf{X}_0 = 10^{2} \mathbf{I}_4.
\end{equation*}
The initial condition $\mathbf{x}_0$ means the pursuer starts $100\mathrm{m}$ behind the evader in the $x$ direction. In the following analysis, we assume each player has perfect information about their own state (i.e., $\mathbf{x}_t^{(p)}$ and $\mathbf{x}_t^{(e)}$) but has partial and asymmetric information about the relative state $\mathbf{x}_t$.

\subsection{Convergence of the Players' Best Response Strategies}\label{sec:optimal_cost}
\begin{figure}[htbp]
    \centering
    \includegraphics[width=.5\textwidth]{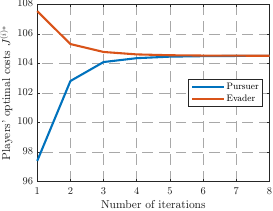}
    \caption{The players' optimal costs converge to a constant value within a few best response iterations.}
    \label{fig:optimal_cost}
\end{figure}

\begin{figure}[htbp]
    \centering
    \includegraphics[width=.5\textwidth]{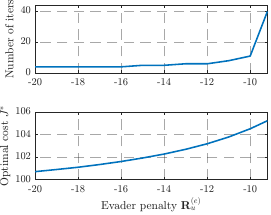}
    \caption{A game's convergence iterations and optimal cost versus evader's input penalty $\mathbf{R}_u^{(e)}$.}
    \label{fig:iteration-opt-cost-r2}
\end{figure}

Figure \ref{fig:optimal_cost} shows the evolution of the optimal costs for both the pursuer and the evader for their augmented belief dynamics at each best response iteration. As the number of iterations increases, each player's optimal cost gradually converges to a constant value. This indicates that both players' best responses no longer improve their individual optimal cost. This point of convergence corresponds to a Nash equilibrium, where neither player obtains a benefit by deviating from their current strategies. Figure \ref{fig:iteration-opt-cost-r2} illustrates that as the penalty magnitude on the evader's input decreases, the number of iterations required for the best response to converge increases. Similarly, the game's optimal value rises with the decreasing penalty magnitude on the evader's input (the pursuer's baseline LQG cost is $97.40$). When the penalty $\mathbf{R}_u^{(e)}$ is reduced to a critical value (the minimum value that prevents the evader from destabilizing the system, which is $-9.2\mathbf{I}_2$ in this case), the number of iterations jumps from $11$ ($\mathbf{R}_u^{(e)} = -10\mathbf{I}_2$) to $40$. In contrast to the big jump in the number of iterations, the game's value increases from $104.52$ to $105.24$. A similar phenomenon was studied in \cite{zhang2021policy} where the $H_\infty$ robustness is preserved in $H_2$ linear control (linear quadratic regulator). In other words, the costs of these problems do not go up to infinity as the eigenvalues of the closed-loop system approach the stability boundary. 

\begin{figure}[htbp]
    \centering
    \includegraphics[width=.5\textwidth]{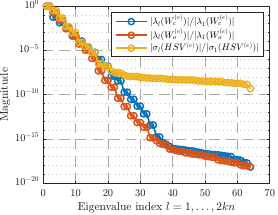}
    \caption{Rapidly decaying eigenvalue magnitudes of the evader's augmented belief dynamics $(k=8,n=4,2kn=64)$, ordered by decreasing magnitude.}
    \label{fig:eigenvalue-decay-rates}
\end{figure}

Figure \ref{fig:eigenvalue-decay-rates} shows the rapid decay rates of Gramian eigenvalues and Hankel singular values' magnitudes in the evader's higher-order belief dynamics (iteration $k=8$). This phenomenon was also observed in \cite{guan2025best} for both players in an open-loop stable system. These rapid decay rates explain why the best response converges after a few iterations: the higher-order belief dynamics become increasingly difficult for the evader to control and observe. Consequently, both players' best responses converge after a few iterations, providing a good approximation of the infinite-dimensional Nash equilibrium. To determine if this phenomenon occurs in general zero-sum LQG dynamic games with partial and asymmetric information, we tested 1000 random open-loop stable systems with higher-dimensional states, inputs, and measurements ($n = 4$ and $m_1 = m_2 = p_1 = p_2 = 2$) than those in \cite{guan2025best}. The results, shown in Table \ref{tab:table_eig}, demonstrate that both players' higher-order belief dynamics are significantly challenging to control and observe in a zero-sum LQG dynamic game with open-loop stable dynamics (each entry in the table shows the percentage of eigenvalues/singular values with magnitudes below a threshold). Therefore, both players' best responses can converge after a few iterations.

\begin{table}[htbp]
    \centering
    \caption{Average proportion of players' controllability Gramian $(\mathbf{W}_c^{(i)})$, observability Gramian $(\mathbf{W}_o^{(i)})$ and Hankel singular values $(\mathbf{HSV}^{(i)})$ less than $\{10^{-4},\ 10^{-8}\}$ of $1000$ random open-loop stable systems.}
    \label{tab:table_eig}
    \renewcommand{\arraystretch}{1.5}
    \begin{tabular}{ccc}
        \hline\hline
         Metrics & $|\lambda_l|,|\sigma_l| \leq 10^{-4}$ & $|\lambda_l|,|\sigma_l| \leq 10^{-8}$ \\
         \hline
        $\mathbf{W}_{c,\mathrm{ave}}^{(p)}$ & $82.33\%$ & $68.84\%$ \\
        $\mathbf{W}_{o,\mathrm{ave}}^{(p)}$ & $81.47\%$ & $69.63\%$ \\
        $\mathbf{HSV}_{\mathrm{ave}}^{(p)}$ & $83.86\%$ & $46.95\%$ \\
        $\mathbf{W}_{c,\mathrm{ave}}^{(e)}$ & $82.82\%$ & $69.89\%$ \\
        $\mathbf{W}_{o,\mathrm{ave}}^{(e)}$ & $80.28\%$ & $69.49\%$ \\
        $\mathbf{HSV}_{\mathrm{ave}}^{(e)}$ & $83.41\%$ & $62.94\%$ \\
        \hline\hline
    \end{tabular}
\end{table}

\subsection{Impacts of the Players' Partial and Asymmetric Beliefs}\label{sec:perf_vs_partial}
\begin{figure}[htbp] 
    \centering
    \begin{subfigure}[t]{.45\textwidth}
        \centering
        \includegraphics[width=\textwidth]{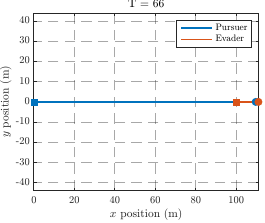}
        \caption{Perfect information (average total cost: $38.47$).}
        \label{fig:perfect_info_traj}
    \end{subfigure}
    \hfill
    \begin{subfigure}[t]{.45\textwidth}
        \centering
        \includegraphics[width=\textwidth]{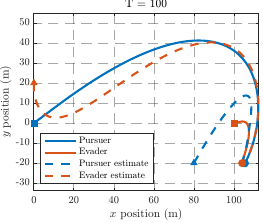}
        \caption{Partial and asymmetric information (average total cost: $104.52$).}
        \label{fig:partial_info_traj}
    \end{subfigure}
    \caption{Mean trajectories (solid lines) of pursuer and evader with their estimates of the opponent's position (dashed lines).}
    \label{fig:perfect_info_vs_partial_info}
\end{figure}

Figure \ref{fig:perfect_info_vs_partial_info} presents the mean trajectories of the pursuer and evader under perfect and partial and asymmetric information scenarios, along with their respective estimates of the opposing player's position. The simulation starts from initial state $(\mathbf{x}_0, \mathbf{z}_0^{(p)}, \mathbf{z}_0^{(e)}) = [(-100, 0, 0, 0), (-80, 20, 0, -25, \dots), (-100, 20, 0, 0, \dots)]$, where higher-order beliefs are identical to the $0$-th order belief. Squares and triangles represent actual and estimated initial positions, respectively, while circles represent final positions. In the perfect information scenario,  as shown in Figure \ref{fig:perfect_info_traj}, the mean trajectories follow a straight line connecting the players' initial positions, as their initial velocities are zero. However, in the partial and asymmetric information scenario, the pursuer's incorrect initial belief about the evader's velocity results in a deviation in the $y$ direction, as shown in Figure \ref{fig:partial_info_traj}. Moreover, the average infinite horizon cost in the partial and asymmetric information scenario is approximately $2.5$ times that of the perfect information scenario, highlighting the impact of partial and asymmetric information on the overall game cost.

\subsection{Impacts of the Players' Belief Orders}\label{sec:belief_order_impact}
\begin{figure}[htbp] 
    \centering
    \begin{subfigure}[t]{.45\textwidth}
        \centering
        \includegraphics[width=\textwidth]{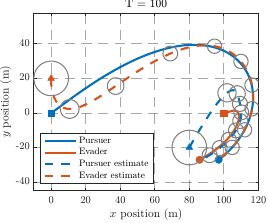}
        \caption{Best response strategies of $1$ iteration (average total cost: $107.54$).}
        \label{fig:mean_traj_low_order_belief}
    \end{subfigure}
    \hfill
    \begin{subfigure}[t]{.45\textwidth}
        \centering
        \includegraphics[width=\textwidth]{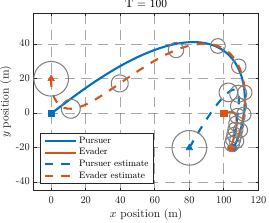}
        \caption{Best response strategies of $8$ iterations (average total cost: $104.52$).}
        \label{fig:mean_traj_higher_order_belief}
    \end{subfigure}
    \caption{Pursuer and evader mean trajectories and estimation uncertainties (ellipses) under best response strategies computed at different iterations.}
    \label{fig:low_order_belief_vs_higher_order_belief}
\end{figure}

Using the best response method, both players can form their beliefs about the mean state and the opposing player's beliefs at different orders through various iterations. The best response strategies based on low- or higher-order beliefs can have different impacts on the game's average infinite horizon cost and their mean trajectories. Figure \ref{fig:low_order_belief_vs_higher_order_belief} illustrates the effects of both players' low- and higher-order beliefs on the game's average infinite horizon cost and their mean trajectories. A notable difference is that the evader travels a longer distance, and the distance between the pursuer and evader is greater when both players use best response strategies based on low-order beliefs compared to higher-order beliefs. Additionally, the game's average infinite horizon cost induced by both players' low-order beliefs is higher than that induced by higher-order beliefs. In this case, the evader prefers both players to have low-order beliefs, while the pursuer prefers higher-order beliefs in terms of optimizing the game's average infinite horizon cost and the distance between the two players.

\subsection{Impacts of the Players' Relative Controllability and Observability}\label{sec:pursuit_evasion_maneuver_observe}
\begin{figure}[htbp]
    \centering
    \includegraphics[width=.5\textwidth]{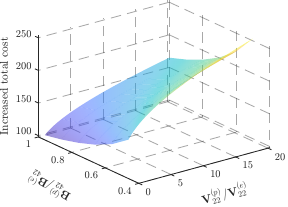}
    \caption{Average infinite horizon cost concerning changes in players' relative controllability $(\mathbf{B}_{42}^{(p)} \in [0.5\Delta t,\ \Delta t])$ and observability $(\mathbf{V}_{22}^{(p)} \in [10^2,\ 45^2])$.}
    \label{fig:total_cost_wo_attack}
\end{figure}

Figure \ref{fig:total_cost_wo_attack} illustrates the game's average infinite horizon cost as it varies with changes in the players' relative controllability (i.e., maneuverability) $(\mathbf{B}_{42}^{(p)} \in [0.5\Delta t,\ \Delta t])$ and observability $(\mathbf{V}_{22}^{(p)} \in [10^2,\ 45^2])$. To ensure the game's upper value remains bounded in all cases, we set the penalty on the evader's dynamic control $\mathbf{R}_u^{(e)} = -100\mathbf{I}_2$, keeping all other parameters the same as in Section \ref{sec:perf_vs_partial}. It is evident that, for the same level of controllability, the game's average infinite horizon cost increases as the pursuer's observation uncertainty rises. Similarly, for the same level of observation uncertainty, the cost rises as the pursuer's controllability decreases. In other words, a pursuer with less controllability and/or noisier observations results in a higher game cost.

\begin{figure}[htbp] 
    \centering
    \begin{subfigure}[t]{.45\textwidth}
        \centering
        \includegraphics[width=\textwidth]{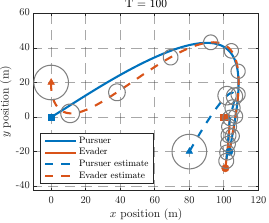}
        \caption{Pursuer and evader's mean trajectories, mean estimates, and estimation uncertainties.}
        \label{fig:less_maneuver_traj}
    \end{subfigure}
    \hfill
    \begin{subfigure}[t]{.45\textwidth}
        \centering
        \includegraphics[width=\textwidth]{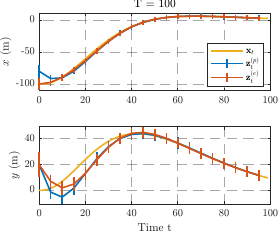}
        \caption{Time-series plots of the actual mean state, pursuer and evader's estimates, and estimation uncertainties.}
        \label{fig:less_maneuver_time_plots}
    \end{subfigure}
    \caption{\textbf{Less controllable pursuer ($\mathbf{B}_{42}^{(p)} = 0.5\Delta t$)}: both players' mean trajectories, $0$-th order estimates, and estimation uncertainties (ellipses) for both players.}
    \label{fig:less_maneuver}
\end{figure}

Figure \ref{fig:less_maneuver} presents the mean trajectories of the pursuer and evader, along with their respective $0$-th order estimates and estimation uncertainties regarding the opposing player's mean state, when the pursuer is less maneuverable $(\mathbf{B}_{42}^{(p)} = 0.5\Delta t)$ than the evader. The average infinite-horizon cost rises to $163.28$ (the pursuer's baseline LQG cost is $147.53$) with $\mathbf{R}_u^{(e)}=-22\mathbf{I}_2$, higher than the baseline value in Section \ref{sec:perf_vs_partial}, reflecting the pursuer's reduced ability to close the relative distance. Additionally, given the same simulation time, the evader travels a longer distance compared to the mean trajectories in Figure \ref{fig:partial_info_traj} due to the pursuer's reduced controllability. Consequently, the pursuer may be unable to capture the evader if there is a specific capture zone (e.g., a safe distance from the base) required for capture.

\begin{figure}[htbp] 
    \centering
    \begin{subfigure}[t]{.45\textwidth}
        \centering
        \includegraphics[width=\textwidth]{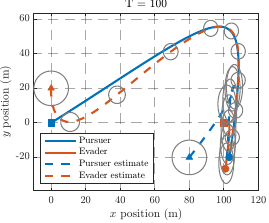}
        \caption{Pursuer and evader's mean trajectories, mean estimates, and estimation uncertainties.}
        \label{fig:noisier_observation_traj}
    \end{subfigure}
    \hfill
    \begin{subfigure}[t]{.45\textwidth}
        \centering
        \includegraphics[width=\textwidth]{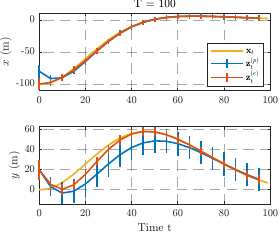}
        \caption{Time-series plots of the actual mean state, pursuer and evader's estimates, and estimation uncertainties.}
        \label{fig:noisier_observation_time_plots}
    \end{subfigure}
    \caption{\textbf{Noisier pursuer observations ($\mathbf{V}_{22}^{(p)} = 37^2$)}: both players' mean trajectories, $0$-th order estimates, and estimation uncertainties (ellipses).}
    \label{fig:noisier_observe}
\end{figure}

Similarly, Figure \ref{fig:noisier_observe} shows the mean trajectories of the pursuer and evader, along with their respective $0$-th order estimates and estimation uncertainties regarding the opposing player's mean state, when the pursuer has noisier observations in the $y$ direction $(\mathbf{V}_{22}^{(p)} = 37^2)$ than the evader. The average infinite horizon cost in this case is $158.74$ (the puruser's baseline LQG cost is $147.43$) with $\mathbf{R}_u^{(e)}=-22\mathbf{I}_2$, which is about $1.5$ times the cost in the partial and asymmetric information scenario in Section \ref{sec:perf_vs_partial}. The pursuer's excursion in the $y$ direction in Figure \ref{fig:noisier_observation_traj} is greater than in Figures \ref{fig:partial_info_traj} and \ref{fig:less_maneuver_traj}. Furthermore, given the same simulation time, the evader travels a significantly longer distance compared to Figure \ref{fig:partial_info_traj} due to the pursuer's increased observation uncertainty in the $y$ direction. 

\begin{figure}[htbp] 
    \centering
    \begin{subfigure}[t]{.45\textwidth}
        \centering
        \includegraphics[width=\textwidth]{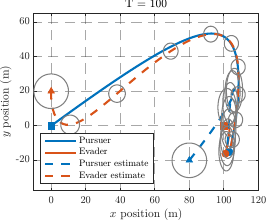}
        \caption{Pursuer and evader's mean trajectories, mean estimates, and estimation uncertainties.}
        \label{fig:better_maneuver_noisier_observe_traj}
    \end{subfigure}
    \hfill
    \begin{subfigure}[t]{.45\textwidth}
        \centering
        \includegraphics[width=\textwidth]{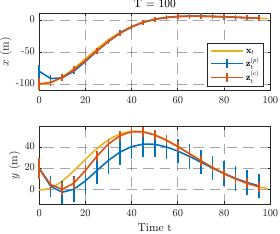}
        \caption{Time-series plots of the actual mean state, pursuer and evader's estimates, and estimation uncertainties.}
        \label{fig:better_maneuver_noisier_observe_time_series}
    \end{subfigure}
    \caption{\textbf{More controllable pursuer with noisier observations ($\mathbf{B}_{42}^{(p)} = 1.5\Delta t, \mathbf{V}_{22}^{(p)} = 37^2$)}: both players' mean trajectories, $0$-th order estimates, and estimation uncertainties (ellipses).}
    \label{fig:better_maneuver_noisier_observe}
\end{figure}

Figure \ref{fig:better_maneuver_noisier_observe}  shows the mean trajectories of the pursuer and evader, along with their respective $0$-th order estimates and estimation uncertainties regarding the opposing player's mean state, when the pursuer is more maneuverable $(\mathbf{B}_{42}^{(p)} = 1.5\Delta t)$ but has noisier observation in the $y$ direction $(\mathbf{V}_{22}^{(p)} = 37^2)$. The average infinite horizon cost in this scenario is $130.21$, which is lower than the cost of $158.74$ when $\mathbf{V}_{22}^{(p)} = 37^2$, but still higher than the cost of $104.52$ in the partial and asymmetric information in Section \ref{sec:perf_vs_partial}. This indicates that while the pursuer's improved controllability can reduce the game cost, the benefits are limited by the increased observation noise in the $y$ direction. This is also evident in the pursuer and evader's mean trajectories in Figure \ref{fig:better_maneuver_noisier_observe_traj}, where the pursuer's excursion in the $y$ direction is still greater than in Figure \ref{fig:partial_info_traj} when both players have the same controllability and observation uncertainty. Moreover, it is possible that the game's upper value becomes unbounded even when the pursuer has a controllability advantage if its observability is sufficiently compromised.

\subsection{Impacts of the Evader's Adversarial Attack}\label{sec:pursuit_evasion_game_attack}

We now investigate the impact of the evader's ECM attack on the pursuer's sensor measurement and the resulting game cost. In Section \ref{sec:two_player_pe_game}, we define the evader's augmented control input $\tilde{\mathbf{u}}_t^{(e)} = [\mathbf{u}_t^{(e)\top},\ \mathbf{s}_t^{(e)\top}]^\top$ that comprises an independent thrust input $\mathbf{u}_t^{(e)}$ affecting the relative kinematics and an ECM input $\mathbf{s}_t^{(e)}$ corrupting only the pursuer's sensor measurement through the feed-through term $\mathbf{D}^{(e)}\mathbf{s}_t^{(e)}$ in \eqref{eq:sensor_ecm}. We set the ECM effectiveness matrix $\mathbf{D}^{(e)} = 4\mathbf{I}_2$, meaning the evader intends to interfere with the pursuer's relative position measurement with gain $4$. The evader's penalty matrices are set to $\mathbf{R}_u^{(e)}=-100\mathbf{I}_2$ and $\mathbf{R}_s^{(e)} = -70\mathbf{I}_2$ that gives the evader a strategic incentive to allocate more effort to sensor attack than to relative kinematics.

\begin{figure}[htbp]
    \centering
    \includegraphics[width=0.5\textwidth]{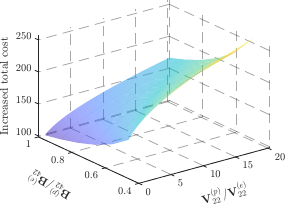}
    \caption{Increased average infinite horizon cost concerning changes in players' relative controllability and observability with both players' convergent best response strategies.}
    \label{fig:total_cost_attack_higher_order_belief}
\end{figure}

\begin{figure}[htbp]
    \centering
    \includegraphics[width=0.5\textwidth]{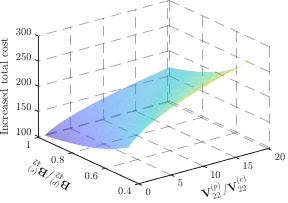}
    \caption{Increased average infinite horizon cost concerning changes in players' relative controllability and observability with both players' best response strategies of $1$ iteration.}
    \label{fig:total_cost_attack_low_order_belief}
\end{figure}

Figures \ref{fig:total_cost_attack_higher_order_belief} and \ref{fig:total_cost_wo_attack} share identical settings for the pursuer's relative controllability $\mathbf{B}_{42}^{(p)}\in[0.5\Delta t,\ \Delta t]$ and observability $\mathbf{V}_{22}^{(p)}=[10^2,\ 45^2]$, differing only in that Figure \ref{fig:total_cost_attack_higher_order_belief} introduces the evader's ECM attack while Figure \ref{fig:total_cost_wo_attack} sets $\mathbf{D}^{(e)}=\mathbf{0}_2$ and $\mathbf{R}_s^{(e)}=\mathbf{0}_2$. When the pursuer's controllability and observability match the evader's ($\mathbf{B}_{42}^{(p)} = \Delta t$, $\mathbf{V}_{22}^{(p)}=10^2$), the game cost increases marginally from $98.03$ without ECM to $98.52$ with ECM, indicating that a pursuer with equal capabilities can largely compensate for degraded measurements through superior estimation. When the pursuer is less controllable and observable than the evader ($\mathbf{B}_{42}^{(p)}=0.5\Delta t$, $\mathbf{V}_{22}^{(p)}=45^2$), the cost rises more sharply from $253.37$ to $256.65$, demonstrating that ECM attacks are most effective against a pursuer already disadvantaged in both controllability and observability. Figure \ref{fig:total_cost_attack_low_order_belief} shows the same ECM scenario but with best response strategies computed at only one iteration rather than the convergent strategies of Figure \ref{fig:total_cost_attack_higher_order_belief}. The minimum and maximum costs in Figure \ref{fig:total_cost_attack_low_order_belief} are $98.57$ and $260.15$ respectively, which is higher than the corresponding convergent values in Figure \ref{fig:total_cost_attack_higher_order_belief} across the entire parameter range. This confirms that convergent best response strategies consistently achieve lower costs than the single-iteration strategies, consistent with the results in Figure \ref{fig:low_order_belief_vs_higher_order_belief}.  

\section{Conclusion}\label{sec:conclusion}
We have formulated and studied a class of zero-sum stochastic LQDGs under partial and asymmetric information. We derived explicit expressions for each player's best response within the class of pure linear dynamic output feedback control strategies. As players iteratively update their best responses, they form increasingly higher-order belief states, leading to infinite-dimensional internal states. However, our numerical results reveal that the game's value converges after only a few iterations, suggesting that higher-order belief states eventually provide no benefit. Thus, an infinite-dimensional Nash equilibrium can be closely approximated by equilibrium strategies based on belief states with limited internal state dimensions. Furthermore, our numerical experiments reveal that factors such as asymmetric beliefs, low-order beliefs, reduced controllability for the pursuer, poor observability, and ECM attacks can significantly increase the game’s cost. To support further research in this area, we provide an open-source implementation of the best response algorithm  \footnote{\url{https://github.com/TSummersLab/partial-asym-info-best-response}}. Although this work focuses on LQDGs, the proposed best response dynamics framework naturally extends to nonlinear high-dimensional pursuit-evasion problems. 

\section*{Appendix: Controllability and Observability Metrics for Higher-Order Belief Dynamics}
This appendix summarizes the controllability and observability analysis from \cite{guan2025best} for the zero-sum LQDG without direct feed-through ($\mathbf{D}^{(1)}=\mathbf{D}^{(2)}=\mathbf{0}$), which the present paper extends to the case $\mathbf{D}^{(i)} \neq \mathbf{0}$.

\subsection*{Controllability and Observability Gramians}
With each best response iteration $k$, the augmented belief dynamics that each player must control and observe become increasingly large-scale. The controllability Gramian $\mathbf{W}_{c,k}^{(i)}$ and observability Gramian $\mathbf{W}_{o,k}^{(i)}$ of player $i$'s augmented belief dynamics at iteration $k$ are defined as the solutions to the discrete-time Lyapunov equations,
\begin{align}
    \mathbf{W}_{c,k}^{(i)} - \mathbf{\bar{A}}_k^{(i)}\mathbf{W}_{c,k}^{(i)}\left( \mathbf{\bar{A}}_k^{(i)} \right)^\top &= \mathbf{\bar{B}}_k^{(i)}\left( \mathbf{\bar{B}}_k^{(i)} \right)^\top, \label{eq:contrl-gramian}\\
    \mathbf{W}_{o,k}^{(i)} - \left( \mathbf{\bar{A}}_k^{(i)} \right)^\top \mathbf{W}_{o,k}^{(i)} \mathbf{\bar{A}}_k^{(i)} &= \left( \mathbf{\bar{C}}_k^{(i)} \right)^\top \mathbf{\bar{C}}_k^{(i)}\label{eq:observe-gramian},
\end{align}
where $\mathbf{W}_{c,k}^{(i)} \succ 0$ and $\mathbf{W}_{o.k}^{(i)} \succ 0$. Eigenvectors of $\mathbf{W}_{c,k}^{(i)}$ and $\mathbf{W}_{o.k}^{(i)}$ associated with small eigenvalues define directions in the augmented state space that are difficult to control and observe respectively.

\subsection*{Hankel Singular Values}
The Hankel singular values of player $i$ at iteration $k$ quantify each state's joint controllability and observability, and are defined as
\begin{equation*}
    \sigma_{j,k}^{(i)} = \sqrt{\lambda_j\left( \mathbf{W}_{c,k}^{(i)}\mathbf{W}_{o,k}^{(i)} \right)}, \quad j=1,2,\dots,n_k^{(i)}n,
\end{equation*}
where $n_k^{(i)}n$ is the dimension of player $i$'s augmented belief state at iteration $k$. Large Hankel singular values correspond to states that are both easy to control and observe; small singular values indicate states that are difficult for the player to exploit strategically.

\subsection*{Cholesky Estimates and Approximation Error Bounds}
To quantify the decay rates analytically, \cite{guan2025best} establishes bounds using a Cholesky factorization of a Cauchy matrix built from the eigenvalues of the augmented belief dynamic matrices $\mathbf{\bar{A}}_k^{(i)}$. Following \cite{guan2025best}, the discrete-time Cholesky factors are defined as
\begin{equation*}
    \delta_{l,k}^{(i)} = \frac{-1}{2 Re\left( \frac{\lambda_{l,k}^{(i)} - 1}{\lambda_{l,k}^{(i)} + 1} \right)} \prod_{j=1}^{l-1} \left\vert \frac{\left(\lambda_{l,k}^{(i)} - \lambda_{j,k}^{(i)}\right)\left(\lambda_{l,k}^{(i)*} + 1\right)}{\left(\lambda_{l,k}^{(i)*}\lambda_{j,k}^{(i)} - 1\right)\left(\lambda_{l,k}^{(i)} + 1\right)} \right\vert^2,
\end{equation*}
where $\lambda_{j,k}^{(i)}\in\mathbf{\lambda}(\mathbf{\bar{A}}_k^{(i)})$ are the eigenvalues of the augmented dynamic matrix, the superscript $*$ denotes complex conjugation, and the eigenvalues are ordered according to a selection rule that yields a descending $\delta_{l,k}^{(i)}$ order.

\textbf{Theorem 1.} Suppose that for each player $i$, the augmented matrices $\mathbf{\bar{A}}_k^i$ are stable and diagonalizable at iteration $k$, the pair $(\mathbf{\bar{A}}_k^{(i)},\ \mathbf{\bar{B}}_k^{(i)})$ is controllable and the pair $(\mathbf{\bar{A}}_k^{(i)},\ \mathbf{\bar{C}}_k^{(i)})$ is observable. Define the low-rank approximation $\hat{\mathbf{W}}_{c_{lm_i},k}^{(i)} = \sum_{j=1}^l \delta_{j,k}^{(i)} \mathbf{Z}_{j,k}^{(i)} \mathbf{Z}_{j,k}^{(i)*}$, where $\mathbf{Z}_{j,k}^{(i)}$ is constructed from the right eigenvectors of $\mathbf{\bar{A}}_k^{(i)}$. If $\delta_{l,k}^{(i)}/\delta_{1,k}^{(i)} < \epsilon$, then $\hat{\mathbf{W}}_{c_{lm_i},k}^{(i)}$ has rank at most $lm_i$ and satisfies 
\begin{equation*}
    \| \mathbf{W}_{c,k}^{(i)} - \hat{\mathbf{W}}_{c_{lm_i},k}^{(i)} \|_\infty \leq \epsilon \delta_{1,k}^{(i)} \left(m_i n_{k}^{(i)}n\right) \left(n_{k}^{(i)}n - l\right) \left( \kappa_\infty \left(\mathbf{X}_k^{(i)}\right) \| \tilde{\mathbf{B}}_k^{(i)} \|_1 \right)^2, 
\end{equation*}
where $\mathbf{X}_k^{(i)}$ denotes the matrix of right eigenvectors of $\mathbf{\bar{A}}_k^{(i)}$ with unit-norm columns, $\kappa_\infty(\cdot)$ is the condition number, $\tilde{\mathbf{B}}_k^{(i)}=\sqrt{2}(\mathbf{\bar{A}}_k^{(i)} + \mathbf{I})\mathbf{\bar{B}}_k^{(i)}$, and $\delta_{1,k}^{(i)}\approx \| \mathbf{W}_{c,k}^{(i)} \|_2$.

Theorem 1 establishes that the controllability Gramians of both players' higher-order belief dynamics can be approximated by low-rank matrices with bounded error, with the bound governed by the ratio $\delta_{l,k}^{(i)}/\delta_{1,k}^{(i)}$. Analogous bounds for the observability Gramian and Hankel singular values follow the same approach. Consequently, Nash equilibrium strategies based on low-order belief states closely approximate the infinite-dimensional equilibrium strategies.

\section*{Funding Sources}
This work was supported by the United States Air Force Office of Scientific Research under Grants FA9550-23-1-0424 and FA2386-24-1-4014, and by the National Science Foundation under Grant ECCS-2047040.

\section*{Acknowledgments}
During the preparation of this manuscript, the author(s) used Claude (Anthropic; versions Sonnet 4.6 and Sonnet 5) to check grammar and improve the clarity and flow of the written content. After using this tool, the author(s) reviewed and edited the content as needed and take full responsibility for the content of the publication.

\bibliography{sample}

@article{bagchi1981linear,
  title={Linear-quadratic stochastic pursuit-evasion games},
  author={Bagchi, Arunabha and Olsder, Geert Jan},
  journal={Applied mathematics and optimization},
  volume={7},
  number={1},
  pages={95--123},
  year={1981},
  publisher={Springer},
  doi={10.1007/BF01442109}
}

@article{molinari1973stable,
  title={The stable regulator problem and its inverse},
  author={Molinari, B},
  journal={IEEE Transactions on Automatic Control},
  volume={18},
  number={5},
  pages={454--459},
  year={1973},
  publisher={IEEE},
  doi={10.1109/TAC.1973.1100364}
}

@article{willems2003least,
  title={Least squares stationary optimal control and the algebraic Riccati equation},
  author={Willems, Jan},
  journal={IEEE Transactions on automatic control},
  volume={16},
  number={6},
  pages={621--634},
  year={2003},
  publisher={IEEE},
  doi={10.1109/TAC.1971.1099831}
}

@article{gupta2014common,
  title={Common information based Markov perfect equilibria for linear-Gaussian games with asymmetric information},
  author={Gupta, Abhishek and Nayyar, Ashutosh and Langbort, C{\'e}dric and Basar, Tamer},
  journal={SIAM Journal on Control and Optimization},
  volume={52},
  number={5},
  pages={3228--3260},
  year={2014},
  publisher={SIAM},
  doi={10.1137/140953514}
}

@article{hambly2023linear,
  title={Linear-quadratic Gaussian Games with Asymmetric Information: Belief Corrections Using the Opponents Actions},
  author={Hambly, Ben and Xu, Renyuan and Yang, Huining},
  journal={arXiv preprint arXiv:2307.15842},
  year={2023},
  doi={10.48550/arXiv.2307.15842}
}

@inbook{isaacsDifferential1999,
  title = {Differential Games: A Mathematical Theory with Applications to Warfare and Pursuit, Control and Optimization},
  shorttitle = {Differential Games},
  author = {Isaacs, Rufus},
  year = {1965},
  publisher = {Courier Corporation},
  chapter={12}
}

@inbook{bacsar1998dynamic,
  title={Dynamic noncooperative game theory},
  author={Ba{\c{s}}ar, Tamer and Olsder, Geert Jan},
  year={1998},
  publisher={SIAM},
  chapter={6},
  doi={10.1137/1.9781611971132}
}

@INPROCEEDINGS{guan2025best,
  author={Guan, Yuxiang and Shames, Iman and Summers, Tyler H.},
  booktitle={2025 European Control Conference (ECC)}, 
  title={Best Response Convergence for Zero-sum Stochastic Dynamic Games with Partial and Asymmetric Information}, 
  year={2025},
  volume={},
  number={},
  pages={2408-2415},
  doi={10.23919/ECC65951.2025.11186954}
}

@inbook{yuksel2024stochastic,
  title={Stochastic teams, games, and control under information constraints},
  author={Y{\"u}ksel, Serdar and Ba{\c{s}}ar, Tamer},
  year={2024},
  publisher={Springer},
  chapter={9,10},
  doi={10.1007/978-3-031-54071-4}
}

@article{harsanyi1967games,
  title={Games with incomplete information played by “Bayesian” players, I--III Part I. The basic model},
  author={Harsanyi, John C},
  journal={Management science},
  volume={14},
  number={3},
  pages={159--182},
  year={1967},
  publisher={INFORMS},
  doi={10.1287/mnsc.14.3.159}
}

@article{behn1968class,
  title={On a class of linear stochastic differential games},
  author={Behn, R and Ho, Yu-Chi},
  journal={IEEE Transactions on Automatic Control},
  volume={13},
  number={3},
  pages={227--240},
  year={1968},
  publisher={IEEE},
  doi={10.1109/TAC.1968.1098898}
}

@article{rhodes1969differential,
  title={Differential games with imperfect state information},
  author={Rhodes, I and Luenberger, D},
  journal={IEEE Transactions on Automatic Control},
  volume={14},
  number={1},
  pages={29--38},
  year={1969},
  publisher={IEEE},
  doi={10.1109/TAC.1969.1099086}
}

@article{rhodes1969stochastic,
  title={Stochastic differential games with constrained state estimators},
  author={Rhodes, I and Luenberger, D},
  journal={IEEE Transactions on Automatic Control},
  volume={14},
  number={5},
  pages={476--481},
  year={1969},
  publisher={IEEE},
  doi={10.1109/TAC.1969.1099281}
}

@article{willman1969formal,
  title={Formal solutions for a class of stochastic pursuit-evasion games},
  author={Willman, Warren},
  journal={IEEE Transactions on Automatic Control},
  volume={14},
  number={5},
  pages={504--509},
  year={1969},
  publisher={IEEE},
  doi={10.1109/TAC.1969.1099249}
}

@article{nayyar2013common,
  title={Common information based Markov perfect equilibria for stochastic games with asymmetric information: Finite games},
  author={Nayyar, Ashutosh and Gupta, Abhishek and Langbort, Cedric and Ba{\c{s}}ar, Tamer},
  journal={IEEE Transactions on Automatic Control},
  volume={59},
  number={3},
  pages={555--570},
  year={2013},
  publisher={IEEE},
  doi={10.1109/TAC.2013.2283743}
}

@article{gupta2016dynamic,
  title={Dynamic games with asymmetric information and resource constrained players with applications to security of cyberphysical systems},
  author={Gupta, Abhishek and Langbort, C{\'e}dric and Ba{\c{s}}ar, Tamer},
  journal={IEEE Transactions on Control of Network Systems},
  volume={4},
  number={1},
  pages={71--81},
  year={2016},
  publisher={IEEE},
  doi={10.1109/TCNS.2016.2584183}
}

@article{pachter2017lqg,
  title={LQG dynamic games with a control-sharing information pattern},
  author={Pachter, Meir},
  journal={Dynamic Games and Applications},
  volume={7},
  pages={289--322},
  year={2017},
  publisher={Springer},
  doi={10.1007/s13235-016-0182-6}
}

@inproceedings{altman2009stochastic,
  title={Stochastic games with one step delay sharing information pattern with application to power control},
  author={Altman, Eitan and Kambley, Vijay and Silva, Alonso},
  booktitle={2009 international conference on game theory for networks},
  pages={124--129},
  year={2009},
  organization={IEEE},
  doi={10.1109/GAMENETS.2009.5137393}
}

@article{tang2023dynamic,
  title={Dynamic games among teams with delayed intra-team information sharing},
  author={Tang, Dengwang and Tavafoghi, Hamidreza and Subramanian, Vijay and Nayyar, Ashutosh and Teneketzis, Demosthenis},
  journal={Dynamic Games and Applications},
  volume={13},
  number={1},
  pages={353--411},
  year={2023},
  publisher={Springer},
  doi={10.1007/s13235-022-00424-4}
}

@article{vasal2021signaling,
  title={Signaling equilibria for dynamic LQG games with asymmetric information},
  author={Vasal, Deepanshu and Anastasopoulos, Achilleas},
  journal={IEEE Transactions on Control of Network Systems},
  volume={8},
  number={3},
  pages={1177--1188},
  year={2021},
  publisher={IEEE},
  doi = {10.1109/TCNS.2021.3059835}
}

@article{ouyang2025approach,
  title={An approach to stochastic dynamic games with asymmetric information and hidden actions},
  author={Ouyang, Yi and Tavafoghi, Hamidreza and Teneketzis, Demosthenis},
  journal={Dynamic Games and Applications},
  volume={15},
  number={1},
  pages={182--215},
  year={2025},
  publisher={Springer},
  doi={10.1007/s13235-024-00558-7}
}

@article{ouyang2016dynamic,
  title={Dynamic games with asymmetric information: Common information based perfect bayesian equilibria and sequential decomposition},
  author={Ouyang, Yi and Tavafoghi, Hamidreza and Teneketzis, Demosthenis},
  journal={IEEE Transactions on Automatic Control},
  volume={62},
  number={1},
  pages={222--237},
  year={2016},
  publisher={IEEE},
  doi={10.1109/TAC.2016.2544936}
}

@inproceedings{sinha2016structured,
  title={Structured perfect Bayesian equilibrium in infinite horizon dynamic games with asymmetric information},
  author={Sinha, Abhinav and Anastasopoulos, Achilleas},
  booktitle={2016 54th Annual Allerton Conference on Communication, Control, and Computing (Allerton)},
  pages={256--263},
  year={2016},
  organization={IEEE},
  doi={10.1109/ALLERTON.2016.7852238}
}

@article{zhang2021policy,
  title={Policy Optimization for $H_2$ Linear Control with $H_\infty$ Robustness Guarantee: Implicit Regularization and Global Convergence},
  author={Zhang, Kaiqing and Hu, Bin and Basar, Tamer},
  journal={SIAM Journal on Control and Optimization},
  volume={59},
  number={6},
  pages={4081--4109},
  year={2021},
  publisher={SIAM},
  doi={10.1137/20M1347942}
}

@article{shima2002time,
  title={Time-varying linear pursuit-evasion game models with bounded controls},
  author={Shima, Tal and Shinar, Josef},
  journal={Journal of Guidance, Control, and Dynamics},
  volume={25},
  number={3},
  pages={425--432},
  year={2002},
  doi = {10.2514/2.4927}
}

@article{shima2011optimal,
  title={Optimal cooperative pursuit and evasion strategies against a homing missile},
  author={Shima, Tal},
  journal={Journal of Guidance, Control, and Dynamics},
  volume={34},
  number={2},
  pages={414--425},
  year={2011},
  doi = {10.2514/1.51765}   
}

@article{weiss2016minimum,
  title={Minimum effort pursuit/evasion guidance with specified miss distance},
  author={Weiss, Martin and Shima, Tal},
  journal={Journal of Guidance, Control, and Dynamics},
  volume={39},
  number={5},
  pages={1069--1079},
  year={2016},
  publisher={American Institute of Aeronautics and Astronautics},
  doi = {10.2514/1.G001623}
}

@article{austin1990game,
  title={Game theory for automated maneuvering during air-to-air combat},
  author={Austin, Fred and Carbone, Giro and Falco, Michael and Hinz, Hans and Lewis, Michael},
  journal={Journal of Guidance, Control, and Dynamics},
  volume={13},
  number={6},
  pages={1143--1149},
  year={1990},
  doi = {10.2514/3.20590}
}

@article{shinar1977analysis,
  title={Analysis of optimal evasive maneuvers based on a linearized two-dimensional kinematic model},
  author={Shinar, Josef and Steinberg, D},
  journal={Journal of Aircraft},
  volume={14},
  number={8},
  pages={795--802},
  year={1977},
  doi = {10.2514/6.1976-1979}
}

@article{pontani2009numerical,
  title={Numerical solution of the three-dimensional orbital pursuit-evasion game},
  author={Pontani, Mauro and Conway, Bruce A},
  journal={Journal of guidance, control, and dynamics},
  volume={32},
  number={2},
  pages={474--487},
  year={2009},
  doi = {10.2514/1.37962}
}

@article{li2020saddle,
  title={Saddle point of orbital pursuit-evasion game under J 2-perturbed dynamics},
  author={Li, Zhen-yu and Zhu, Hai and Yang, Zhen and Luo, Ya-zhong},
  journal={Journal of Guidance, Control, and Dynamics},
  volume={43},
  number={9},
  pages={1733--1739},
  year={2020},
  publisher={American Institute of Aeronautics and Astronautics},
  doi = {10.2514/1.G004459}
}

@article{huang2021dynamic,
  title={A dynamic game framework for rational and persistent robot deception with an application to deceptive pursuit-evasion},
  author={Huang, Linan and Zhu, Quanyan},
  journal={IEEE Transactions on Automation Science and Engineering},
  volume={19},
  number={4},
  pages={2918--2932},
  year={2021},
  publisher={IEEE},
  doi={10.1109/TASE.2021.3097286}
}

@inproceedings{shishika2021partial,
  title={Partial information target defense game},
  author={Shishika, Daigo and Maity, Dipankar and Dorothy, Michael},
  booktitle={2021 IEEE International Conference on Robotics and Automation (ICRA)},
  pages={8111--8117},
  year={2021},
  organization={IEEE},
  doi={10.1109/ICRA48506.2021.9561995}
}

@inproceedings{cavalieri2014incomplete,
  title={Incomplete information pursuit-evasion games with uncertain relative dynamics},
  author={Cavalieri, Kurt A and Satak, Neha and Hurtado, John E},
  booktitle={AIAA Guidance, Navigation, and Control Conference},
  pages={0971},
  year={2014},
  doi={10.2514/6.2014-0971}
}

@article{scassellati2002theory,
  author  = {Scassellati, Brian},
  title   = {Theory of Mind for a Humanoid Robot},
  journal = {Autonomous Robots},
  volume  = {12},
  number  = {1},
  pages   = {13--24},
  year    = {2002},
  doi = {https://doi.org/10.1023/A:1013298507114}
}

@inproceedings{swarup2003linear,
  title={Linear-Quadratic-Gaussian differential games with different information patterns},
  author={Swarup, Ashitosh and Speyer, Jason L},
  booktitle={42nd IEEE International Conference on Decision and Control (IEEE Cat. No. 03CH37475)},
  volume={4},
  pages={4146--4151},
  year={2003},
  organization={IEEE},
  doi={10.1109/CDC.2003.1271799}
}

@article{chong1971stochastic,
  title={On the stochastic control of linear systems with different information sets},
  author={Chong, Chee-Yee and Athans, Michael},
  journal={IEEE Transactions on Automatic Control},
  volume={16},
  number={5},
  pages={423--430},
  year={1971},
  publisher={IEEE},
  doi={10.1109/TAC.1971.1099810}
}

@book{antoulas2005approximation,
  title={Approximation of large-scale dynamical systems},
  author={Antoulas, Athanasios C},
  year={2005},
  publisher={SIAM},
  doi={10.1137/1.9780898718713}
}

\end{document}